\documentclass[preprint,a4paper,10pt,dvipsnames]{article}
\usepackage[T1]{fontenc}
\usepackage[utf8]{inputenc}
\usepackage[english]{babel}
\usepackage{amsmath,amsthm,enumerate}
\usepackage{amssymb}
\usepackage{vmargin}
\usepackage{url}
\usepackage{stmaryrd}
\usepackage{tikz}
\usepackage{tikz-network}
\usetikzlibrary{decorations.pathreplacing}
\usepackage{subfig}
\usepackage{authblk}

\tikzstyle{v}=[circle,inner sep=0, minimum size =6 pt, line width = 1pt, draw=black, fill=black, text= white]

\newcommand{\qedclaim}{\hfill $\diamond$ \medskip}

\begin{document}

\title{The Maker-Breaker Largest Connected Subgraph Game\thanks{This work has been supported by the European Research Council (ERC) consolidator grant No.~725978 SYSTEMATICGRAPH, the National Research Agency (ANR) project P-GASE (ANR-21-CE48-0001-01), the STIC-AmSud project GALOP, the PHC Xu Guangqi project DESPROGES, and the UCA$^\textsc{jedi}$ Investments in the Future project managed by the National Research Agency (ANR-15-IDEX-01).}}

\author[1]{Julien Bensmail}
\author[1]{Foivos Fioravantes}
\author[2]{Fionn Mc~Inerney}
\author[1]{Nicolas Nisse}
\author[3,4]{Nacim Oijid}

\affil[1]{Universit\'e C\^ote d'Azur, CNRS, Inria, I3S, Biot, France}
\affil[2]{CISPA Helmholtz Center for Information Security, Saarbr\"{u}cken, Germany}
\affil[3]{\'Ecole Normale Sup\'erieure de Lyon, 69364 Lyon Cedex 07, France }
\affil[4]{Univ. Lyon, Universit\'e Lyon 1, LIRIS UMR CNRS 5205, F-69621, Lyon, France}

\maketitle

\begin{abstract}
Given a graph $G$ and $k \in \mathbb{N}$, we introduce the following game played in $G$. Each round, Alice colours an uncoloured vertex of $G$ red, and then Bob colours one blue (if any remain). Once every vertex is coloured, Alice wins if there is a connected red component of order at least $k$, and otherwise, Bob wins. This is a Maker-Breaker version of the Largest Connected Subgraph game introduced in [Bensmail et al. The Largest Connected Subgraph Game. {\it Algorithmica}, 84(9):2533--2555, 2022]. We want to compute $c_g(G)$, which is the maximum $k$ such that Alice wins in $G$, regardless of Bob's strategy.

Given a graph $G$ and $k\in \mathbb{N}$, we prove that deciding whether $c_g(G)\geq k$ is PSPACE-complete, even if $G$ is a bipartite, split, or planar graph. To better understand the Largest Connected Subgraph game, we then focus on {\it A-perfect} graphs, which are the graphs $G$ for which $c_g(G)=\lceil|V(G)|/2\rceil$, {\it i.e.}, those in which Alice can ensure that the red subgraph is connected. We give sufficient conditions, in terms of the minimum and maximum degrees or the number of edges, for a graph to be A-perfect. Also, we show that, for any $d \geq 4$, there are arbitrarily large A-perfect $d$-regular graphs, but no cubic graph with order at least $18$ is A-perfect. Lastly, we show that $c_g(G)$ is computable in linear time when $G$ is a $P_4$-sparse graph (a superclass of cographs). 
\end{abstract}

\newtheorem{theorem}{Theorem}[section]
\newtheorem{lemma}[theorem]{Lemma}
\newtheorem{proposition}[theorem]{Proposition}
\newtheorem{conjecture}[theorem]{Conjecture}
\newtheorem{observation}[theorem]{Observation}
\newtheorem{claim}{Claim}
\newtheorem{corollary}[theorem]{Corollary}
\newtheorem{question}[theorem]{Question}


\section{Introduction}

{\it Maker-Breaker} games have been vastly studied since the introduction of some of the famous Maker-Breaker games like {\it Hex}, introduced independently by Hein and Nash in the 1940s~\cite{G59}, and the {\it Shannon switching game} of Shannon from the 1950s~\cite{G61}. Maker-Breaker games drew more attention after the 1973 paper of Erdős and Selfridge on {\it positional} games~\cite{ES73}, a superclass of Maker-Breaker games (see~\cite{HKSS14} for more on positional games). They now form a major domain in combinatorial game theory, and more largely, theoretical computer science.

In Maker-Breaker games, there are a set of elements $X$, and a family of winning sets $\mathcal{F}$, which is a family of subsets of $X$. The two players, {\it Maker} and {\it Breaker}, alternate selecting previously unselected elements of $X$. Maker wins if she selects every element of a winning set in $\mathcal{F}$, while Breaker wins if he prevents this, {\it i.e.}, by selecting at least one element of each winning set in $\mathcal{F}$. One of the first major results for Maker-Breaker games is the Erdős-Selfridge Theorem~\cite{ES73} from 1973, which gives sufficient conditions for Breaker to win. In 1978, Schaefer proved that determining the winner of a Maker-Breaker game is PSPACE-complete, even if each of the winning sets in $\mathcal{F}$ has size at most $11$~\cite{Schaefer1978}. This was not improved upon until 2021, when Rahman and Watson proved that the same holds even if each of the winning sets in $\mathcal{F}$ has size at most $6$ (or exactly $6$)~\cite{rahman}. The latter results significantly impacted complexity theory, as the Maker-Breaker game (equivalently, POS CNF) is a classic problem to reduce from in order to prove PSPACE-hardness, and the size of the winning sets often has implications on the properties of the hard instances of the problem being reduced to. In terms of applications, the study of such games has also often led to {\it a priori} unrelated results. For example, the general method of the proof of the Erdős-Selfridge Theorem~\cite{ES73} was the first instance of the method of conditional expectations, the first technique to efficiently derandomise randomised algorithms~\cite{HKSS14}. Also, one of the cornerstone papers on Ramsey theory is actually on positional games~\cite{HJ63}.

Apart from general results for Maker-Breaker games, many individual such games have been considered. Some of the more notable ones were introduced by Chv\'{a}tal and Erdős, and are played on the complete graph $K_n$: the {\it Hamiltonicity game}, the {\it Connectivity game}, and the {\it Clique game}~\cite{CE78}. In each of these games, $X$ consists of the edges of $K_n$, while $\mathcal{F}$ consists of all Hamiltonian cycles for the former, all spanning trees for the second, and all cliques of a given size for the latter. They notably also introduced {\it biased} Maker-Breaker games, in which Breaker selects multiple elements of $X$ on each of his turns, and the goal is to determine the least number he may select, while still guaranteeing him a win~\cite{CE78}.

Along the same lines, in this paper, we introduce the following Maker-Breaker game, which is a natural game that has, surprisingly, not been considered in the literature to date. In the {\it Maker-Breaker Largest Connected Subgraph game} played on a given graph $G$, there is a positive integer $k$ given as an input, and $X$ consists of the vertices of $G$, while $\mathcal{F}$ consists of all the connected subgraphs of order at least $k$ in $G$. In particular, for a given graph $G$, we are interested in the parameter $c_g(G)$, which is the largest integer $k$ such that Maker has a winning strategy in the Maker-Breaker Largest Connected Subgraph game in $G$.  

Another motivation for introducing this game is that it is a Maker-Breaker version of the {\it Largest Connected Subgraph game} of Bensmail et al.~\cite{paper1}. The Largest Connected Subgraph game is played by two players, Alice and Bob\footnote{Although referring to the two players of a Maker-Breaker game as Maker and Breaker is standard, to preserve the connection of this game with the Largest Connected Subgraph game, we instead refer to them as Alice and Bob.}, through successive rounds played on a graph $G$. Initially, every vertex of $G$ is uncoloured. Each round, Alice colours an uncoloured vertex of $G$ red, and then, Bob colours an uncoloured vertex blue (if any remain).
The game ends once all the vertices of $G$ have been coloured, resulting in a {\it red subgraph} of $G$ (induced by the vertices coloured red) and a {\it blue subgraph} of $G$ (induced by the vertices coloured blue).
Alice's (Bob's, resp.) score is the order (number of vertices) of the largest connected component of the red subgraph (blue subgraph, resp.). If the players have different scores, then the player with the largest score wins. Otherwise, the game ends in a draw. In~\cite{paper1}, it was proved, through standard strategy-stealing arguments, that Alice always has a strategy to ensure at least a draw, and thus, that Bob can never win if Alice plays optimally. It was also proved that determining the outcome of the game ({\it i.e.}, determining the winner of the game or if the game ends in a draw) on a given graph is PSPACE-complete, even when restricted to bipartite graphs of diameter~5, but polynomial-time solvable for paths, cycles, and cographs. This game was novel as it is natural, and there is also a rich background on these types of games, {\it i.e.}, it is a connection game (see~\cite{c05} for a book on these games), since the players seek to make connected structures, and a {\it scoring} game (see~\cite{LNS17} for a survey on these games), as the players' scores determine the winner. Thus, our game is also a connection game, and is closely related to a scoring game.

We now further corroborate the introduction of our game in relation to the Largest Connected Subgraph game. In the Largest Connected Subgraph game, for some subgraphs of certain graphs, Bob prefers to limit Alice's score in them, rather than increase his score in them. In particular, this can be true in disconnected graphs. Bob limiting Alice's score in these subgraphs is equivalent to playing the Maker-Breaker Largest Connected Subgraph game in them. Another motivation is to understand the properties of graphs in which Alice can ensure a single connected red component at the end of the Largest Connected Subgraph game (especially since Alice wins in these graphs if they have odd order). Thus, we study {\it A-perfect graphs}, which are the graphs $G$ for which $c_g(G)=\lceil |V(G)|/2\rceil$, {\it i.e.}, the graphs in which Alice can ensure a single connected red component at the end of the (Maker-Breaker) Largest Connected Subgraph game. Thus, some results in this paper can be directly applied to the Largest Connected Subgraph game, but as many of our techniques are based on the structural properties of the graphs considered, our methods could be applied to other Maker-Breaker games.

Indeed, other games that are relatively close to ours can also be found in the literature, such as biased positional games~\cite{biased1} (where Maker aims to maximise the order of the largest connected component of his graph through selecting edges of a complete graph), or variants of them played on vertices similar to as in our game~\cite{biased3,biased2}. Specifically, there is a close relation between our game and the {\it Maker-Breaker Domination game}~\cite{DGPR20} (see Section~\ref{section:complexity}). In the Maker-Breaker Domination game, $X$ consists of the vertices of the graph, while $\mathcal{F}$ consists of all the dominating sets of the graph~\cite{DGPR20}. Deciding the winner of the Maker-Breaker Domination game is PSPACE-complete, even for bipartite graphs and split graphs, but can be determined in polynomial time for cographs and trees~\cite{DGPR20}. Determining how fast Maker can win in this game has also been studied~\cite{GIK19}, as well as a variant where $\mathcal{F}$ consists of all the total dominating sets of the graph~\cite{FM22,GHIK20}. However, our game is closer to a connected variant of the Maker-Breaker Domination game, where $\mathcal{F}$ consists of all connected dominating sets instead, but, to the best of our knowledge, this variant has yet to be studied (although there is a connected version of the Domination game~\cite{BFS19}). Indeed, it is easy to see that (see Lemma~\ref{lemma connected dominating set}), given a graph $G$, if Maker wins the connected variant of the Maker-Breaker Domination game, then $G$ is A-perfect. However, the other direction of this statement does not necessarily hold, as can be seen in the following example. Let $G$ be a star of odd order, {\it i.e.}, a graph consisting of one universal vertex $u$ adjacent to an even number of degree-$1$ vertices. Let $G'$ be the graph obtained from $G$ by subdividing exactly one of the edges of $G$. Then, $G'$ is A-perfect since Alice will first colour the vertex $u$ that is universal in $G$, and since $G'$ is of even order, then Bob plays last, and so, Alice can never be forced to colour the degree-$1$ vertex that is at distance~$2$ from $u$. It is clear that Maker does not win the connected variant of the Maker-Breaker Domination game in $G'$. However, we do wonder if this relation holds for graphs of odd order. That is, does it hold that, given a graph $G$ of odd order, Maker wins the connected variant of the Maker-Breaker Domination game if and only if $G$ is A-perfect? Such a result would be a nice connection between our game and this variant of the Maker-Breaker Domination game.    

\paragraph{\bf Our contributions} In Section~\ref{section:preliminaries}, the main terminology and early observations, to be used throughout, are introduced. In Section~\ref{section:complexity}, we show that, given a graph $G$ and an integer $k\geq 1$, deciding whether $c_g(G)\geq k$ is PSPACE-complete, even if $G$ is restricted to be in the class of bipartite graphs of diameter~$4$, split graphs, or planar graphs. In Section~\ref{section:bounds}, we give two sufficient conditions (in terms of the minimum and maximum degrees, and the number of edges) for a graph to be A-perfect. We also prove that arbitrarily large A-perfect $d$-regular graphs exist if and only if $d\geq 4$. In particular, any A-perfect cubic ({\it i.e.}, 3-regular) graph has order at most 16. In Section~\ref{section:sparse}, we show that $c_g$ can be determined in linear time for $(q,q-4)$-graphs, a superclass of cographs. We conclude in Section~\ref{section:conclusion} with a discussion including perspectives for further work on the topic.


\section{Preliminaries}\label{section:preliminaries}

\subsection{Graph theory terminology and notation}

Throughout this paper, all the graphs we consider are undirected and simple.
For a graph $G$, we denote by $V(G)$ its set of vertices, and by $E(G)$ its set of edges.
For a vertex $v$ of $G$, we denote by $N_G(v)$ its {\it neighbourhood}, which is the set of vertices that are adjacent to $v$ in $G$.
The {\it closed neighbourhood} of $v$, denoted by $N_G[v]$, is the set $\{v\} \cup N_G(v)$.
These two notions of neighbourhood extend to sets $S$ of vertices of $G$,
with $N_G(S)$ referring to the subset of vertices of $V(G) \setminus S$ that have a neighbour in $S$,
and $N_G[S]$ referring to the set $S \cup N_G(S)$.
For a vertex $v$ of $G$, we denote by $d_G(v)$ its {\it degree}, with $d_G(v)=|N_G(v)|$. When the graph $G$ is clear from the context, we will drop the subscript in the parameters $N_G$ and $d_G$,
and write them as $N$ and $d$ instead.
The parameters $\delta(G)$ and $\Delta(G)$ refer to the {\it minimum degree} and {\it maximum degree}, respectively, of a vertex in $G$.

For a set $S$ of vertices or edges of $G$, we denote by $G-S$ the {\it subgraph} of $G$ resulting from the deletion of the elements in $S$. If $S=\{x\}$, we will write $G-x$ instead of $G-S$.
Similarly, we denote by $G+S$ (or $G+x$, for short, if $S=\{x\}$) the {\it supergraph} of $G$ obtained by adding the elements (vertices or edges) of $S$.
We denote by $G[S]$ the subgraph of $G$ {\it induced} by the elements in $S$.

A vertex $u$ of $G$ {\it dominates} another vertex $v$ if $v\in N_G(u)$.
We say that $u$ is {\it universal} if $N_G[u]=V(G)$.
A set $S$ of vertices of $G$ is {\it dominating} if $N_G[S]=V(G)$.
This dominating set $S$ is {\it connected} if $G[S]$ is connected. The {\it distance} between two vertices $u$ and $v$ of $G$ is the length of a shortest path from $u$ to $v$. The {\it diameter} of $G$ is the maximum distance between two vertices of $G$. For any standard notion or terminology on graphs not defined in this work, see~\cite{D12}.

\subsection{Additional terminology for the game}

Due to the close ties of our game with the Largest Connected Subgraph game, we will refer to Maker as Alice, and Breaker as Bob. To also make the distinction of who selected which vertices easier, we will say that Alice colours a vertex red when she selects a vertex, and that Bob colours a vertex blue when he selects a vertex. Furthermore, we will refer to the score of Alice as the largest connected red component in the graph at the end of the game. Thus, for a given graph $G$, $c_g(G)$ is the maximum score for which Alice has a strategy ensuring at least this score in $G$. 

A {\it strategy} for a player $P$ is a function $\mathcal{S}$ taking all the previous moves of both players (and the order of these moves, hence, the history of the game) as an input, and outputting the next move for player $P$. Given a graph $G$, an {\it optimal strategy for Alice} is one that ensures her a score of at least $c_g(G)$, while an {\it optimal strategy for Bob} is one that ensures Alice's score is at most $c_g(G)$. Since our game is a {\it parity game}~\cite{K02}, optimal strategies can actually be determined from just the current configuration of coloured vertices, rather than also knowing the order these vertices were coloured in. Thus, for our game, there can also be an equivalent (in terms of optimality) second definition of a strategy for a player $P$, which is a function $\mathcal{S}$ that takes the current configuration of coloured vertices and outputs the next move for player $P$. We will interchangeably use both definitions, depending on which one suits us best at the time.

Throughout this paper, several of our proofs rely on the fact that Alice or Bob can reach a certain game configuration ({\it i.e.}, have a certain set of vertices coloured with their colour) early on. In such cases, to lighten the exposition,
we will sometimes allow ourselves to expose only the most important moves of the strategies that Alice or Bob must make in some rounds of the game. In particular, the reader should keep in mind that, in each of the strategies we describe, 1) if Alice or Bob cannot colour a given vertex in a given round because that vertex is already coloured, then they must colour any other uncoloured vertex instead, and
2) if no vertex to colour for Alice or Bob in a given round is specified, then they must colour any uncoloured vertex. Note that arbitrary moves will never hinder the accomplishments of the main thread of a strategy for that player.

\subsection{General results and observations}

First, we show that the parameter $c_g$ is closed under taking subgraphs, and that when playing the game on a disconnected graph, Alice should focus on the connected component which is the most favourable for her.

\begin{lemma} \label{lemma subgraph}
If $H$ is a (not necessarily proper) subgraph of a graph $G$ with connected components $H_1,\dots,H_k$, then $c_g(H)=\max\left\{c_g(H_1),\dots,c_g(H_k)\right\} \le c_g(G)$.
\end{lemma}

\begin{proof}
First, we show that the parameter $c_g$ is closed under taking subgraphs. We give a strategy for Alice ensuring her a score of at least $c_g(H)$ in $G$. Alice first plays in $H$ according to an optimal strategy $\mathcal{S}$ in $H$. Then, whenever Bob plays in $H$, Alice responds in $H$ according to $\mathcal{S}$, and if this is not possible (the vertex to be coloured by $\mathcal{S}$ is already coloured or there are no uncoloured vertices in $H$) or Bob plays in $G$, then Alice colours an arbitrary uncoloured vertex in $G$. In particular, whenever Alice is forced to colour an arbitrary uncoloured vertex in $H$, she ignores the fact that vertex is coloured when considering her strategy $\mathcal{S}$ in $H$ in the future. The result follows since Alice will colour at least all the vertices in $H$ that she would colour by $\mathcal{S}$, ensuring her a score of at least $c_g(H)$ in $G$ since $\mathcal{S}$ is optimal in $H$.

Then, we show that when playing the game on a disconnected graph, Alice should focus on the connected component which is the most favourable for her. That is, if $G$ is a graph with connected components $H_1,\dots,H_k$, then $c_g(G) = \max\left\{c_g(H_1),\dots,c_g(H_k)\right\}.$ The lower bound follows from the paragraph above, and the upper bound holds since the $k$ connected components are pairwise disconnected, so Bob can just respond in the same connected component that Alice just played in each time (when this is not possible, he colours an arbitrary uncoloured vertex in $G$, which can only be beneficial to him).
\end{proof}





As can be seen in the following corollary, Lemma~\ref{lemma subgraph} implies that a disconnected A-perfect graph must be of even order and have exactly two connected components, one of which consists of a single vertex.

\begin{corollary}\label{cor:A-perfect}
If $G$ is an A-perfect graph that is not connected, then $G$ is of even order and consists of two connected components, one of which consists of a single vertex.    
\end{corollary}

\begin{proof}
Since $G$ is not connected, it consists of connected components $H_1,\ldots,H_k$ for some $k\geq 2$. If $k\geq 3$ and/or $\min\left\{|V(H_1)|,\dots,|V(H_k)|\right\}\geq 2$, then $\max\left\{|V(H_1)|,\dots,|V(H_k)|\right\}\leq |V(G)|-2$, and so, $\left\lceil\frac{|V(G)|-2}{2}\right\rceil\geq \max\left\{c_g(H_1),\dots,c_g(H_k)\right\}=c_g(G)$ by Lemma~\ref{lemma subgraph}, {\it i.e.}, $G$ is not A-perfect. Hence, $k=2$ and $\min\left\{|V(H_1)|,|V(H_2)|\right\}=1$. If $G$ is of odd order, then w.l.o.g., let $|V(H_1)|\geq 2$ be even and $|V(H_2)|=1$. Then, by Lemma~\ref{lemma subgraph}, $c_g(G)=\max\left\{c_g(H_1),c_g(H_2)\right\}\leq \frac{|V(H_1)|}{2}<\left\lceil\frac{|V(G)|}{2}\right\rceil$.
\end{proof}

As will be seen later on, Alice can exploit different types of strategies to achieve the best possible score for her. One such strategy, that is particularly relevant in sufficiently dense graphs, is through colouring the vertices of a connected dominating set.

\begin{lemma}\label{lemma connected dominating set}
For a graph $G$, if, at any point in the game, Alice has coloured all the vertices of a connected dominating set of $G$, then her score will be $\left \lceil \frac{|V(G)|}{2} \right \rceil$.
\end{lemma}

\begin{proof}
Assume Alice has coloured the vertices of a connected dominating set $S$ at some point in the game. By the connectivity property of $S$, there must be, once the game ends, a connected red component containing the vertices of $S$.
Also, by the dominating property of $S$, all the vertices of $G$ not in $S$ have at least one neighbour in $S$.
This implies that the red subgraph must be connected, and thus, that Alice achieves a score of $\left \lceil \frac{|V(G)|}{2} \right \rceil$.
\end{proof}


\section{Computational complexity}\label{section:complexity}

Recall that in~\cite{paper1}, the Largest Connected Subgraph game was shown to be PSPACE-complete, even when restricted to bipartite graphs of diameter~$5$. In this section, using a similar reduction scheme, we prove that the Maker-Breaker Largest Connected Subgraph game is also PSPACE-complete, that is, given a graph $G$ and an integer $k
\geq 1$, deciding whether $c_g(G)\geq k$ is PSPACE-complete.
In fact, we prove that our game is PSPACE-complete, even when restricted to bipartite graphs of diameter~$4$,
split graphs, or planar graphs.

Similarly as in~\cite{paper1},
we establish some of our PSPACE-completeness results via reductions from POS CNF, 
a game for which deciding whether Alice or Bob has a winning strategy was shown to be PSPACE-complete in \cite{Schaefer1978}. 
This game is a 2-player game where the input $(X,\phi)$ consists of a set of variables $X = \{x_1, \dots , x_n \}$, 
and of a formula $\phi$ in conjunctive normal form (CNF) made up of clauses $C_1, \dots , C_m$ each containing variables of $X$ in their positive forms.
Each round, the first player, Alice, sets a variable of $\phi$ (that is not yet set) to true, before the second player, Bob, sets a variable of $\phi$ (that is not yet set) to false. 
Once all the variables of $X$ have been assigned a truth value, Alice wins if $\phi$ is true, and Bob wins if $\phi$ is false.

As mentioned earlier, our game is closely related to the Maker-Breaker Domination game~\cite{DGPR20}. For this reason, the PSPACE-hardness reduction we give for proving the upcoming Theorem~\ref{theorem:pspace-bipartite} (Corollary~\ref{corollary:pspace-split}, resp.) is very similar to (the same as, resp.) the one given in~\cite{DGPR20}, for proving similar complexity results. However, we include the full proofs for completeness. It should be noted that these are rather standard reductions, but we later give a more clever one to prove our game is PSPACE-complete in planar graphs. 

Before we start, note first that when given a graph $G$ and an integer $k\geq 1$, the problem of deciding whether $c_g(G)\geq k$ is in PSPACE since, in the game, there are $\lceil|V(G)|/2\rceil$ rounds and the number of possible moves for a player in a round is at most $|V(G)|$. Thus, in the upcoming proofs, we focus on proving the PSPACE-hardness of the game.

\begin{theorem}\label{theorem:pspace-bipartite}
Given a graph $G$ and an integer $k\geq 1$, it is PSPACE-complete to decide whether $c_g(G)\geq k$, even when $G$ is restricted to be in the class of bipartite graphs of diameter~$4$.
\end{theorem}

\begin{proof}
We prove the PSPACE-hardness via a reduction from POS CNF. Let $(X, \phi)$ be an instance of POS CNF. 
Set $X = \{x_1, \dots, x_n\}$ and $\phi = C_1 \wedge \dots \wedge C_m$. 
By adding a dummy variable in $X$ if needed, we can suppose $n$ is even.

Consider the graph $G$ constructed as follows. For every variable $x_i \in X$, we add a vertex $v_i$ to $G$. 
For every clause $C_j$ of $\phi$, we add two vertices $C_j^1$ and $C_j^2$ to $G$. 
For every variable $x_i \in X$ and clause $C_j$ of $\phi$,
we add the edges $v_iC_j^1$ and $v_iC_j^2$ to $G$ if $x_i$ appears in $C_j$.
Finally, we add two vertices $u_1$ and $u_2$ to $G$, that we make adjacent to all of the $v_i$'s. Note that the resulting $G$, which is constructed in polynomial time, 
is bipartite and has diameter at most~$4$.

Set $k = |V(G)|/2$, and note that $|V(G)|$ is even. We will show that Alice wins in $(X,\phi)$ if and only if $c_g(G) \geq k$. Let us assume first that Alice has a winning strategy in $(X,\phi)$.
We give a strategy for Alice that ensures a score of at least $k$ when playing the Maker-Breaker Largest Connected Subgraph game in $G$. In the first round, Alice colours the vertex $v_i$ that corresponds to the variable $x_i \in X$ she would have set to true in the first round of her winning strategy in $(X,\phi)$. From the second round on, in each round, if the last vertex coloured by Bob is

\begin{itemize}
    \item some $v_i$, then Alice colours the vertex $v_j$ corresponding to the variable $x_j$ she would set to true in response to Bob setting $x_i$ to false in her winning strategy in $(X, \phi)$; 
    	
	\item $u_1$ ($u_2$, resp.), then Alice colours $u_2$ ($u_1$, resp.);     	    	
    	
    \item some $C_j^1$ ($C_j^2$, resp.), then Alice colours $C_j^2$ ($C_j^1$, resp.).
\end{itemize}

Whenever Alice cannot follow the strategy above, she colours an arbitrary uncoloured vertex. By Alice's strategy, once the game in $G$ ends, exactly one vertex in every pair $\{C_j^1,C_j^2\}$ is red, exactly one vertex in $\{u_1,u_2\}$ is red, and the $v_i$'s corresponding to the $x_i$'s she would have set to true in her winning strategy for $(X,\phi)$ are also red. 
Because Alice wins in $(X,\phi)$ with that strategy, every vertex $C_j^{\ell}$ of $G$ coloured red must be adjacent to at least one vertex $v_k$ coloured red corresponding to a variable she would have set to true when playing in $(X,\phi)$. Since all the $v_i$'s are dominated by $u_1$ and $u_2$, and one of these two vertices is red, we deduce that the red subgraph must contain only one connected component. Thus, Alice achieves a score of $k$ and $c_g(G) \ge k$.

Assume now that Bob has a winning strategy in $(X,\phi)$. We give a strategy for Bob that ensures that Alice's score is strictly less than $k$ when playing the Maker-Breaker Largest Connected Subgraph game in $G$. In each round, if the last vertex coloured by Alice is 

\begin{itemize}
	\item some $v_i$, then Bob colours the vertex $v_j$ corresponding to the variable $x_j$ he would set to false in response to Alice setting $x_i$ to true in his winning strategy in $(X,\phi)$;
	\item $u_1$ ($u_2$, resp.), then Bob colours $u_2$ ($u_1$, resp.);
	\item some $C_j^1$ ($C_j^2$, resp.), then Bob colours $C_j^2$ ($C_j^1$, resp.).
\end{itemize}

Note that Bob can follow this strategy from start to end, as $n$ is even. By Bob's strategy, once the game in $G$ ends, exactly one vertex in every pair $\{C_j^1, C_j^2\}$ is red.
Also, since Bob coloured all the $v_i$'s corresponding to $x_i$'s he would set to false when following a winning strategy in $(X,\phi)$, there exists a $C_{q}$ that is not satisfied in $(X,\phi)$, meaning its variables were all set to false by Bob. In $G$, this translates to exactly one of $C_q^1$ or $C_q^2$ being red while all of their neighbours (the $v_i$'s corresponding to the $x_i$'s that $C_q$ contains), are blue. Thus, the red subgraph contains at least two connected components, and hence, Alice achieves a score of less than~$k$, and $c_g(G) < k$.
\end{proof}

\begin{corollary}\label{corollary:pspace-split}
Given a graph $G$ and an integer $k\geq 1$, it is PSPACE-complete to decide whether $c_g(G)\geq k$, even when $G$ is restricted to be in the class of split graphs.
\end{corollary}

\begin{proof}
The proof is similar to that of Theorem~\ref{theorem:pspace-bipartite}, with the slight difference being in the construction of $G$. Here, neither of the vertices $u_1$ and $u_2$ are added, while all the possible edges between the $v_i$'s are added so that they form a clique, thus making $G$ a split graph. The same strategies for Alice and Bob (omitting $u_1$ and $u_2$) from the proof of Theorem~\ref{theorem:pspace-bipartite} remain applicable by the same arguments, and the result follows.
\end{proof}

From the proofs of Theorem~\ref{theorem:pspace-bipartite} and Corollary~\ref{corollary:pspace-split}, it follows that deciding if a bipartite graph with diameter $4$ (resp., a split graph) is A-perfect is PSPACE-complete. To prove that the game is PSPACE-complete, even when restricted to planar graphs, we need a different reduction. This time, we establish the result by a reduction from {\it Planar Generalised Hex}, which was proved to be PSPACE-complete~\cite{Reisch1981}. Planar Generalised Hex is played on a planar graph $G$, in which a particular {\it outside pair} $\{s,t\}$ of vertices, {\it i.e.}, $st \not \in E(G)$ and $G+st$ is planar, is set. Initially, $s$ and $t$ are red.
Then, in successive rounds, the first player, Alice, colours an uncoloured vertex red, before the second player, Bob, then colours an uncoloured vertex blue. The game ends once all the vertices of $G$ have been coloured. If the red subgraph contains a path joining $s$ and $t$, then Alice wins. Otherwise, Bob wins. We can now prove our last result in this section.

\begin{figure}
    \centering
    \scalebox{0.85}{
	\begin{tikzpicture}
		\Vertex[x=0,y=0, fontscale = 2, size = 3.4, style = {color = gray!40}, fontcolor=white]{4}
		
		\draw (0,0) node(){$H$};
		\draw (-1.7,0) node[v](t){}  node[right=0.2] {$s$};
		\draw (1.7,0) node[v](s){}  node[left=0.2] {$t$};
		\draw (-2.7,0) node[v](t0){}  node[below right=0.1] {$s_0^1$};
		\draw (2.7,0) node[v](s0){}  node[below left=0.1] {$t_0^1$};
		
		\Edge[](t)(t0)
		\Edge[](s)(s0) 
		
		\draw (-3.7,-.8) node[v](t1){}  node[below=0.2] {$s_1^1$} node[above=.5]{\vdots};
		\draw (-3.7,.8) node[v](tn){}  node[above=0.2] {$s_{n+4}^1$};

		\draw (3.7,-.8) node[v](s1){}  node[below=0.2] {$t_1^1$} node[above=.5]{\vdots};;
		\draw (3.7,.8) node[v](sn){}  node[above=0.2] {$t_{n+4}^1$};

		\Edge[](t1)(t0)
		\Edge[](tn)(t0)
		\Edge[](s1)(s0)
		\Edge[](sn)(s0)
		
		\draw (-1.7,1.7) node[v](t0){}  node[below left=0.1] {$s_0^2$};
		\draw (1.7,1.7) node[v](s0){}  node[below right=0.1] {$t_0^2$};
		
		\Edge[](t)(t0)
		\Edge[](s)(s0) 
		
		\draw (-0.9,2.7) node[v](t1){}  node[right=0.2] {$s_{n+4}^2$} node[left = .4]{\dots};
		\draw (-2.5,2.7) node[v](tn){}  node[left=0.2] {$s_1^2$};

		\draw (0.9,2.7) node[v](s1){}  node[left = 0.2] {$t_1^2$}node[right = .4]{\dots};
		\draw (2.5,2.7) node[v](sn){}  node[right=0.2] {$t_{n+4}^2$};

		\Edge[](t1)(t0)
		\Edge[](tn)(t0)
		\Edge[](s1)(s0)
		\Edge[](sn)(s0)
		 
		\draw (-1.7,-1.7) node[v](t0){}  node[above left=0.1] {$s_0^3$};
		\draw (1.7,-1.7) node[v](s0){}  node[above right=0.1] {$t_0^3$};
		
		 \Edge[](t)(t0)
		 \Edge[](s)(s0) 
		
		\draw (-0.9,-2.7) node[v](t1){}  node[right=0.2] {$s^3_{n+4}$}node[left = .4]{\dots};
		\draw (-2.5,-2.7) node[v](tn){}  node[left=0.2] {$s^3_1$};

		\draw (0.9,-2.7) node[v](s1){}  node[left=0.2] {$t^3_1$}node[right = .4]{\dots};
		\draw (2.5,-2.7) node[v](sn){}  node[right=0.2] {$t^3_{n+4}$};

		\Edge[](t1)(t0)
		\Edge[](tn)(t0)
		\Edge[](s1)(s0)
		\Edge[](sn)(s0)
	\end{tikzpicture}
}    
    \caption{Illustration of the construction in the proof of Theorem~\ref{theorem:pspace-planar}.}
    \label{fig:planar reduction}
\end{figure}
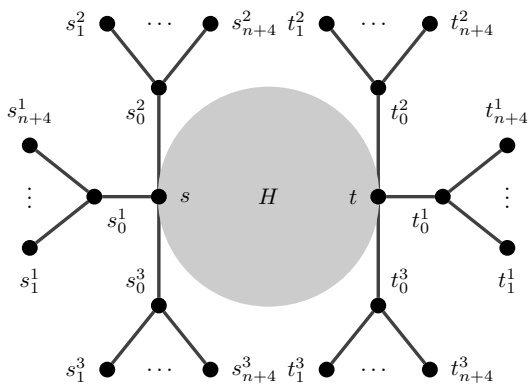

\begin{theorem}\label{theorem:pspace-planar}
Given a graph $G$ and an integer $k\geq 1$, it is PSPACE-complete to decide whether $c_g(G)\geq k$, even when $G$ is restricted to be in the class of planar graphs.
\end{theorem}

\begin{proof}
We prove the PSPACE-hardness via a reduction from Planar Generalised Hex. Let $(H,s,t)$ be an instance of Planar Generalised Hex such that $H$ is the planar graph with the outside pair $\{s,t\}$, that the game is being played on. Set $n = |V(H)|$. By adding a degree-$1$ vertex (a leaf) in $H$ if needed, we can suppose $n$ is even, as this will not change the outcome of $(H,s,t)$. Let $G$ be the graph constructed as follows (see Figure~\ref{fig:planar reduction}). Start from $G$ being the graph $H$. Then, add three vertices $s^1_0, s^2_0,s^3_0$ and make each of them adjacent to $s$, and add another three vertices $t^1_0, t^2_0,t^3_0$, and make each of those adjacent to $t$. Finally, to each of these six vertices we have just added, attach $n+4$ new degree-$1$ vertices, so that a total of $6(n+4)$ degree-$1$ vertices (leaves) are added to $G$. The construction is achieved in polynomial time, and since $H$ is planar, $G$ is too.

Set $k = n + 5$. We will show that Alice wins in $(H,s,t)$ if and only if $c_g(G)\geq k$. Let us assume first that Alice has a winning strategy in $(H,s,t)$. We give a strategy for Alice that ensures a score of at least $k$ when playing the Maker-Breaker Largest Connected Subgraph game in $G$. In the first round, Alice colours~$s$. In the second round, Alice colours $s^1_0$ if possible, and if not, then she colours $s^2_0$. From the third round on,

	\begin{itemize}
		\item if Alice can colour a vertex in $\{s^1_0, s^2_0, s^3_0\}$ in the third round, then she does. If so, then, in each of the next rounds, if possible, Alice colours an uncoloured neighbour of an $s^i_0$ she coloured earlier. At the end of the game, the red subgraph will contain a connected component of order at least $3 + \left \lceil \frac{2(n+4) - 3}{2} \right \rceil = n + 6$, and thus, Alice will have a score of at least~$k$;
		
		\item otherwise, Bob coloured two vertices in $\{s^1_0, s^2_0, s^3_0\}$ in the first two rounds. Then, Alice colours $t$ in the third round, and she then colours one of $t^1_0$ and $t^2_0$ in the fourth round. At this point, for the same reasons as earlier, if Bob has not coloured two vertices in $\{t^1_0, t^2_0, t^3_0\}$ by the end of the fourth round,
		then Alice can colour a vertex in that set in the fifth round, and, as above, guarantee herself a score of at least~$k$.
	\end{itemize}
	
	Thus, we can suppose that, after four rounds, w.l.o.g., $s$, $t$, $s^1_0$, and $t^1_0$ are red, while $s^2_0$, $s^3_0$, $t^2_0$, and $t^3_0$ are blue. From here, Alice's strategy continues as follows. In the fifth round, Alice colours, in $G$, the vertex of $H$ she would have coloured in the first round of her winning strategy in $(H,s,t)$. From the sixth round on, in each round, if the last vertex coloured by Bob in $G$ is
	
	\begin{itemize}
		\item some vertex $u\in V(H)$, 
		then Alice colours, in $G$, the vertex of $H$ she would have coloured in her winning strategy in $(H,s,t)$, as an answer to Bob colouring $u$;
		\item a leaf adjacent to some $s^i_0$ or $t^i_0$, then Alice colours another uncoloured leaf adjacent to the same vertex.
	\end{itemize}

Whenever Alice cannot follow the strategy above, she colours an arbitrary uncoloured vertex. By this strategy, at the end of the game in $G$, $s$ and $t$ are red, and all the vertices that Alice would have coloured through her winning strategy in $(H,s,t)$ are also red. Moreover, $s^1_0$ and $t^1_0$ are red, and, for each of them, she coloured half of their adjacent leaves. Thus, the red subgraph contains a connected component of order at least $n+8$. Thus, Alice achieves a score of at least~$k$, and $c_g(G) \ge k$.

Assume now that Bob has a winning strategy in $(H,s,t)$. We give a strategy for Bob that ensures that Alice's score is strictly less than $k$ when playing the Maker-Breaker Largest Connected Subgraph game in $G$. In each round, if the last vertex coloured by Alice is 

\begin{itemize}
    \item in $\{s,s^1_0,s^2_0,s^3_0\}$, then Bob colours a vertex in $\{s,s^1_0,s^2_0,s^3_0\}$;
    \item in $\{t,t^1_0,t^2_0,t^3_0\}$, then Bob colours a vertex in $\{t,t^1_0,t^2_0,t^3_0\}$;
    \item a vertex $u$ of $H-\{s,t\}$, then Bob colours the vertex of $H$ he would have coloured by his winning strategy in $(H,s,t)$, as an answer to Alice colouring $u$;
    \item a leaf adjacent to some $s^i_0$ or $t^i_0$, then Bob colours another uncoloured leaf adjacent to the same vertex.
\end{itemize}

Note that Bob always answers to one of Alice's moves by colouring a vertex in a set with even size since $n$ is even. Thus, Bob can follow this strategy from start to end. At the end of the game in $G$, the largest connected component of the red subgraph cannot contain both $s$ and $t$, as the moves made by Alice and Bob correspond exactly to the moves that would have been made if they had played in $(H,s,t)$. Moreover, there cannot be two $s^i_0$'s belonging to the same connected red component, as, by the strategy above, Bob must have coloured $s$ in this case. The same holds for the $t^i_0$'s. Also, for any of the $s^i_0$'s and $t^i_0$'s, by Bob's strategy above, Alice can have coloured at most half of the leaves adjacent to it. Thus, because Alice coloured at most half of the vertices in $H-\{s,t\}$, the largest connected red component in $G$ must have order at most $\frac{n-2}{2} + 2 + \frac{n+4}{2} = n+3$. 
Thus, Alice achieves a score of less than $k$, and $c_g(G) < k$.
\end{proof}


\section{A-perfect graphs and regular graphs}\label{section:bounds}

In this section, we study graphs $G$ for which $c_g(G)$ equals one of the next bounds:

\begin{lemma}\label{lemma bounds}
For every graph $G$, $\left \lfloor \frac{\Delta(G)}{2} \right \rfloor + 1 \le c_g(G) \le \left \lceil \frac{|V(G)|}{2} \right \rceil.$
\end{lemma}

\begin{proof}
The right-hand side of the inequality follows from the fact that Alice always colours exactly $\left \lceil \frac{|V(G)|}{2} \right \rceil$ vertices. We now give a strategy for Alice that ensures a score of at least $\left \lfloor \frac{\Delta(G)}{2} \right \rfloor + 1$, to prove the left-hand side of the inequality. In the first round, Alice colours a vertex $v$ with degree $\Delta(G)$. Then, in each of the next rounds, if possible, Alice colours an uncoloured neighbour of $v$. Once the game ends, by the strategy above, Alice must have coloured $v$ and at least half of its neighbours, and the result follows.
\end{proof}

Both bounds in Lemma~\ref{lemma bounds} can be reached for arbitrarily large graphs. For the upper bound, recall that a graph $G$ is {\it A-perfect} if $c_g(G)=\lceil |V(G)|/2 \rceil$. For example, there exist arbitrarily large connected graphs that are A-perfect, since every graph with a universal vertex is A-perfect. Regarding the lower bound, the graph $G$ that is the disjoint union of $m$ copies of the complete graph $K_{d+1}$ (for any $d\in \mathbb{N}$) is $d$-regular, and $c_g(G)=\left \lfloor \frac{d}{2}\right \rfloor + 1$, while $G$ gets more and more distant from being A-perfect as $m$ increases. In particular, since Alice wins the Largest Connected Subgraph game in A-perfect graphs of odd order, we seek conditions for graphs to be A-perfect. We give two sufficient conditions, one based on the minimum and maximum degrees, and the other based on the number of edges.

\begin{theorem} \label{thm large degree}
If $G$ is a connected graph with $\Delta(G) + \delta(G) \geq |V(G)|$, then $G$ is A-perfect.
\end{theorem}

\begin{proof}
We give a strategy for Alice ensuring that, at the end of the game, the red subgraph is connected, which implies that $G$ is A-perfect. Let $u$ be any vertex of degree $\Delta(G)$. In the first round, Alice colours $u$. For every $i\geq 1$, let $C_i$ be the connected component of red vertices at the end of the $i^{th}$ round (we will show that the red vertices always induce a connected subgraph, and so, $C_i$ is well-defined). Let $R_i = V(G) \setminus N[C_i]$, {\it i.e.}, $R_i$ is the set of (non-red) vertices not dominated by a red vertex at the end of the $i^{th}$ round, and let $R^U_i$ be the subset of uncoloured vertices in $R_i$ at the end of the $i^{th}$ round. Note that $C_1=\{u\}$ is connected and that 
$$|R^U_1| \leq |R_1| = |V(G)| - \left|N[C_1]\right| = |V(G)| - \Delta(G) - 1 \leq \delta(G) - 1.$$

Let us show by induction on $i\geq 1$ that, at the end of the $i^{th}$ round, $C_i$ is connected and either $R^U_i= \emptyset$ (in which case we are done) or $|R^U_i| \leq \delta(G) - i$. By the above paragraph, the induction hypothesis holds for $i=1$. Let $i\geq 1$ and let us assume that the induction hypothesis holds for $i$. We show it still holds for $i+1$.

If $R^U_i=\emptyset$, then $C_i$ is a connected red dominating set of the subgraph of $G$ induced by the vertices of $C_i$ and the remaining uncoloured vertices of $G$. From now on, Alice may colour any uncoloured vertex, and the induction hypothesis clearly holds for $i+1$. In particular, the set of red vertices induces a connected subgraph at the end of the game, proving the result.

Otherwise, let $v \in R^U_i$. Since $v$ has at least $\delta(G)$ neighbours (none of which are red since $N[R_i] \cap C_i = \emptyset$) and Bob has coloured $i$ vertices, $v$ has at least $\delta(G)-i$ uncoloured neighbours, and $\delta(G)-i>0$ since $R^U_i\neq \emptyset$ and $|R^U_i|\leq \delta(G)-i$. Moreover, $|R^U_i \setminus \{v\}|<\delta(G)-i$, so $v$ has at least one uncoloured neighbour $w$ not in $R_i$, which implies that $w \in N(R_i) = N(C_i)$. In the $(i+1)^{th}$ round, Alice colours $w$. Then, $C_{i+1}=C_i \cup \{w\}$ is clearly connected, and $R^U_{i+1} \subseteq R_{i+1} \subseteq R_i\setminus \{v\}$ (since $v \in N(C_{i+1})$), and hence, $|R^U_{i+1}| \leq |R_{i+1}| \leq |R_i|-1\leq \delta(G)-(i+1)$.
\end{proof}

We note that the bound in the statement of Theorem~\ref{thm large degree} is sharp, in the sense that there exists a graph $G$ with $\Delta(G)+\delta(G) = |V(G)|-1$ that is not A-perfect. For example, consider the graph $G$ consisting of two complete graphs on $d \geq 3$ vertices joined by a single edge $e$. Then, $\Delta(G)=d$, $\delta(G)=d-1$, $|V(G)|=2d$, and thus, $\Delta(G)+\delta(G)=2d-1=|V(G)|-1$. However, Bob can guarantee that Alice achieves a score of about $|V(G)|/4$, by colouring an uncoloured vertex incident to $e$ in the first round, and then, in each subsequent round, colouring an uncoloured vertex in the same clique that Alice just coloured a vertex in. Thus, $G$ is not A-perfect.

The next result shows that if $G$ has sufficiently many edges, then $G$ is A-perfect.

\begin{theorem}\label{thm:dense} 
If $G$ is a connected graph with $|E(G)|>\frac{(|V(G)|-2)(|V(G)|-3)}{2}+2$, then $G$ is A-perfect.
\end{theorem}

\begin{proof}
Set $n=|V(G)|$, $m=|E(G)|$, and $$x=\frac{(n-2)(n-3)}{2}+2=\frac{n^2-5n+10}{2}.$$
Note first that $\Delta(G)\geq n-4$. Indeed, if we had $\Delta(G)\leq n-5$, then we would deduce that $m \leq \frac{n(n-5)}{2}<x$, which contradicts that $m>x$. Furthermore, if $\Delta(G)=n-4$, then $\delta(G)\geq 7$. Indeed, if there is a degree-$6$ vertex, then we have a contradiction since $$m\leq \frac{(n-1)(n-4)+6}{2}=\frac{n^2-5n+10}{2}=x.$$ Thus, if $\Delta(G)=n-4$, then $G$ is A-perfect by Theorem~\ref{thm large degree} since $\delta(G)\geq 7$. Lastly, if $\Delta(G)=n-1$, then $\delta(G)\geq 1$, and thus, $G$ is A-perfect by Theorem~\ref{thm large degree}. Hence, in what follows, we assume that $n-3\leq \Delta(G)\leq n-2$. We give a strategy for Alice that allows her to colour the vertices of a connected dominating set of $G$ within the first four rounds, and so, by Lemma~\ref{lemma connected dominating set}, $G$ is A-perfect. We treat the two possible values for $\Delta(G)$ independently.

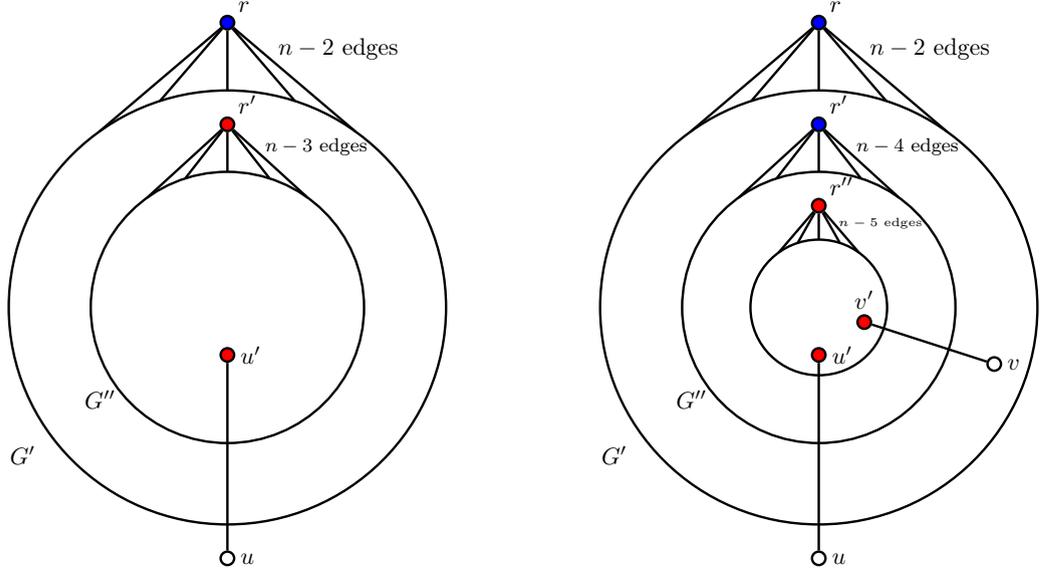
\begin{figure}[!t]
\centering
\subfloat[{The state of the game in $G$ after Alice's first two turns in Case 1.(a), where $G''=G'[N(r')]$.}]{
\scalebox{0.9}{
	\begin{tikzpicture}[inner sep=0.7mm]
		\draw[-,line width=1pt] (0,0) circle (3.2) ;
		
		\node[draw, circle, line width=1pt, fill=blue](r) at (2.5*360/10:3.2cm+1cm)[label=above right: $r$]{};
		
		\node(G') at (6*360/10:3.2cm+0.5cm)[]{$G'$};
		\node(dr) at (1.86*360/10:3.2cm+0.95cm)[]{$n-2$ edges};
		\node[draw, circle, line width=1pt, fill=red](u') at (7.5*360/10:3.2cm-2.5cm)[label=right: $u'$]{};
		\node[draw, circle, line width=1pt, fill=white](u) at (7.5*360/10:3.2cm+0.5cm)[label=right: $u$]{};
		
		\draw[-, line width=1pt]  (r) -- (2*360/10:3.2cm);
		\draw[-, line width=1pt]  (r) -- (3*360/10:3.2cm);
		\draw[-, line width=1pt]  (r) -- (3.5*360/10:3.2cm);
		\draw[-, line width=1pt]  (r) -- (2.5*360/10:3.2cm);
		\draw[-, line width=1pt]  (r) -- (1.5*360/10:3.2cm);
		\draw[-, line width=1pt]  (u) -- (u');

		\draw[-,line width=1pt] (0,0) circle (2) ;
		\node[draw, circle, line width=1pt, fill=red](r') at (2.5*360/10:1.5cm+1.2cm)[label=above right: $r'$]{};

		\draw[-, line width=1pt]  (r') -- (2*360/10:2cm);
		\draw[-, line width=1pt]  (r') -- (3*360/10:2cm);
		\draw[-, line width=1pt]  (r') -- (3.5*360/10:2cm);
		\draw[-, line width=1pt]  (r') -- (2.5*360/10:2cm);
		\draw[-, line width=1pt]  (r') -- (1.5*360/10:2cm);
		\node(G'') at (6*360/10:2cm+0.3cm)[]{$G''$};
		\node(dr') at (1.7*360/10:2cm+0.7cm)[]{\footnotesize{$n-3$ edges}};
	\end{tikzpicture}
}
}
\hspace{40pt}
\subfloat[The state of the game in $G$ after Alice's first three turns in Case 1.(b)i.]{
\scalebox{0.9}{
\begin{tikzpicture}[inner sep=0.7mm]
	\draw[-,line width=1pt] (0,0) circle (3.2) ;
	
	\node[draw, circle, line width=1pt, fill=blue](r) at (2.5*360/10:3.2cm+1cm)[label=above right: $r$]{};
	
	\node(G') at (6*360/10:3.2cm+0.5cm)[]{$G'$};
	\node(dr) at (1.86*360/10:3.2cm+0.95cm)[]{$n-2$ edges};
	\node[draw, circle, line width=1pt, fill=red](u') at (7.5*360/10:3.2cm-2.5cm)[label=right: $u'$]{};
	\node[draw, circle, line width=1pt, fill=white](u) at (7.5*360/10:3.2cm+0.5cm)[label=right: $u$]{};
	
	\draw[-, line width=1pt]  (r) -- (2*360/10:3.2cm);
	\draw[-, line width=1pt]  (r) -- (3*360/10:3.2cm);
	\draw[-, line width=1pt]  (r) -- (3.5*360/10:3.2cm);
	\draw[-, line width=1pt]  (r) -- (2.5*360/10:3.2cm);
	\draw[-, line width=1pt]  (r) -- (1.5*360/10:3.2cm);
	\draw[-, line width=1pt]  (u) -- (u');
	
	\draw[-,line width=1pt] (0,0) circle (2) ;
	\node[draw, circle, line width=1pt, fill=blue](r') at (2.5*360/10:1.5cm+1.2cm)[label=above right: $r'$]{};
	
	\draw[-, line width=1pt]  (r') -- (2*360/10:2cm);
	\draw[-, line width=1pt]  (r') -- (3*360/10:2cm);
	\draw[-, line width=1pt]  (r') -- (3.5*360/10:2cm);
	\draw[-, line width=1pt]  (r') -- (2.5*360/10:2cm);
	\draw[-, line width=1pt]  (r') -- (1.5*360/10:2cm);
	\node(G'') at (6*360/10:2cm+0.3cm)[]{$G''$};
	\node(dr') at (1.7*360/10:2cm+0.7cm)[]{\footnotesize{$n-4$ edges}};
	
	\draw[-,line width=1pt] (0,0) circle (1) ;
	\node[draw, circle, line width=1pt, fill=red](r'') at (2.5*360/10:1.5cm)[label=above right: $r''$]{};

	\node[draw, circle, line width=1pt, fill=white](v) at (9.5*360/10:3.2cm-0.5cm)[label=right: $v$]{};
	\node[draw, circle, line width=1pt, fill=red](v') at (9.5*360/10:3.2cm-2.5cm)[label=above: $v'$]{};

	\draw[-, line width=1pt]  (r'') -- (2*360/10:1cm);
	\draw[-, line width=1pt]  (r'') -- (3*360/10:1cm);
	\draw[-, line width=1pt]  (r'') -- (3.5*360/10:1cm);
	\draw[-, line width=1pt]  (r'') -- (2.5*360/10:1cm);
	\draw[-, line width=1pt]  (r'') -- (1.5*360/10:1cm);
	\draw[-, line width=1pt]  (v') -- (v);
	
	\node(dr'') at (1.5*360/10:1+1.5cm)[]{\tiny{$n-5$ edges}};
\end{tikzpicture}
}
}
\caption{Cases 1.(a) and 1.(b)i. in the proof of Theorem~\ref{thm:dense}.}
\label{figure:dense1}
\end{figure}

\begin{enumerate}
	\item $\Delta(G)=n-2$. 
	
	Let $r\in V(G)$ be such that $d(r)=n-2$, and let $G'=G[N(r)]$. Then, $|V(G')|=n-2$. Since $d(r)=n-2$, there is exactly one additional vertex $u\in V(G)\setminus V(G')$ ($u\neq r$). If $d(u)\geq 2$, then Alice colours $r$ in the first round, and then, in the second round, she colours a neighbour of $u$ (this is possible since $d(u)\geq 2$), and these vertices form a connected dominating set of $G$. Thus, we may assume that $d(u)=1$, and let $N(u)=\{u'\}$. We have that $\Delta(G')\geq n-4$. Indeed, if $\Delta(G')\leq n-5$, then we have a contradiction since $$m\leq \frac{(n-2)(n-5)}{2}+n-2+1=\frac{n^2-5n+8}{2}<x.$$ 
  
	We distinguish the following subcases: 
	\begin{enumerate}
		\item $\Delta(G')=n-3$ (see Figure~\ref{figure:dense1}(a) for an illustration). 
		
		Let $r'\in V(G')$ be such that $d_{G'}(r')=n-3$. Alice's strategy is as follows.
		She starts by colouring $u'$. Now, if Bob does not colour $r$, 
		then Alice continues by colouring $r$, at which point she has coloured the vertices of the connected dominating set $\{u',r\}$ of $G$.
		So, we may assume that Bob colours $r$ in the first round. In the second round, Alice colours $r'$. Observe that $\{u',r'\}$ also forms a connected dominating set of $G$ since $d_{G'}(r')=n-3$, and thus, $u'r'\in E(G)$.  
    
		\item $\Delta(G')=n-4$. 
		
		Let $r'\in V(G')$ be such that $d_{G'}(r')=n-4$, and let $G''=G'[N(r')]$. We distinguish cases according to whether $u'\in V(G'')$ or not. 
    
		\begin{enumerate}
			\item $u'\in V(G'')$. 
			
			Since $d_{G'}(r')=n-4$, there is exactly one additional vertex $v\in V(G')\setminus V(G'')$ ($v\neq r'$). Note that $d_{G'}(v)\geq 1$ because if $d_{G'}(v)=0$, {\it i.e.}, $N(v)=r$, then we have a contradiction since $$m\leq \frac{(n-3)(n-4)}{2}+n-2+1=\frac{n^2-5n+10}{2}=x.$$
			If $d_{G'}(v)\geq 2$, then Alice's strategy is as follows. 
			She starts by colouring $u'$. As before, Bob is forced to colour $r$ in the first round. In the second round, Alice colours $r'$. In the third round, if $u'\notin N(v)$, then Alice colours a neighbour $v'\in V(G')$ (this is possible since $d_{G'}(v)\geq 2$).  After three rounds, Alice's vertices form a connected dominating set of $G$. 
        
			Assume now that $d_{G'}(v)=1$, and let $N_{G'}(v)=\{v'\}$ (see Figure~\ref{figure:dense1}(b) for an illustration). Then, $\Delta(G'')=n-5$. Indeed, if $\Delta(G'')\leq n-6$ (and so, $n\geq 6$), then we have a contradiction since $$m\leq \frac{(n-4)(n-6)}{2}+n-4+n-2+1+1=\frac{n^2-6n+16}{2}\leq x.$$ 
			Let $r''\in V(G'')$ be such that $d_{G''}(r'')=n-5$, and observe that $v'\in N(r'')$. Alice's strategy is as follows. She starts by colouring $u'$, forcing Bob to colour $r$. Then, she colours $v'$ forcing Bob to colour $r'$ (similarly to earlier, if Bob does not colour $r'$, then Alice colours $r'$, and thus, has coloured the vertices of a connected dominating set of $G$). Finally, Alice colours $r''$. Observe that the vertices $u'$, $v'$, and $r''$ form a connected dominating set of $G$.   
        
			\item $u'\notin V(G'')$. 
			
			Observe that $u'$ is the only vertex of $G'$ that is not a neighbour of $r'$, and that $d(u')\geq 3$. Indeed, if $d(u')\leq 2$, then we have a contradiction since $$m\leq \frac{(n-3)(n-4)}{2}+n-2+1=\frac{n^2-5n+10}{2}=x.$$
			Thus, there is at least one edge $u'u''$ with $u''\in V(G'')$. If $d(u')\geq 4$, then Alice's strategy is as follows. She starts by colouring $u'$, forcing Bob to colour $r$.	Then, she colours $r'$, and, in the third round, she colours one of the remaining uncoloured neighbours of $u'$ in $G''$ (which exists since $d(u')\geq 4$). These three vertices form a connected dominating set of $G$. 
			
			Otherwise, $d(u')=3$, and, as in Case 1.(b)i, there exists $r''\in V(G'')$ such that $d_{G''}(r'')=n-5$. Alice's strategy is as follows. She starts by colouring $u'$, forcing Bob to colour $r$.	Then, she colours $u''$, forcing Bob to colour $r'$. 
			Finally, Alice colours $r''$. Note that $u'$, $u''$, and $r''$ form a connected dominating set of $G$. 
		\end{enumerate} 
	\end{enumerate}

	\item $\Delta(G)=n-3$. 
	
    Observe that $G$ cannot contain two vertices $u,v$ such that $d(u)+d(v)\leq 5$. Indeed, if there are two such vertices, then we have a contradiction since $m\leq \frac{(n-2)(n-3)+5}{2}$, but this is not an integer since $(n-2)(n-3)+5$ is odd, and thus, $$m\leq \frac{(n-2)(n-3)+5-1}{2}=\frac{n^2-5n+10}{2}=x.$$ 
    Let $r$ be a vertex of $G$ such that $d(r)=n-3$, and let $G'=G[N(r)]$. Since $d(r)=n-3$, there are exactly two additional vertices $u,v\in V(G)\setminus V(G')$ ($u,v\neq r$). We distinguish cases according to the degrees of $u$ and $v$, and note that $d(u)+d(v)\geq 6$. In what follows, when we say that Alice colours a vertex if needed, it means that if it is not necessary (in the sense that such a vertex has already been coloured), then she either colours the vertex she is supposed to colour in the next round, or she colours an arbitrary uncoloured vertex in that round. 
    
    \begin{enumerate}
    
    \item $d(u),d(v)\geq 3$. 
        
    Alice's strategy is as follows. She starts by colouring $r$. In the second round, she colours an uncoloured vertex in $N(v)$ in $G'$ (this is possible since $d(v)\geq 3$). In the third round, if needed, {\it i.e.}, if Alice has not coloured a vertex in $N(u)$ yet, Alice colours an uncoloured vertex in $N(u)$ in $G'$ if possible, and if not, then $uv\in E(G)$ and Bob coloured $N(u)\setminus \{v\}$ in the first two rounds, and so, she colours $v$. Then, by the end of the third round, Alice has coloured $r$, at least one vertex in $N(v)$ in $G'$, and at least one vertex in $N(u)$, and these vertices form a connected dominating set of $G$.
    
    \item $d(u)=2$ and $d(v)\geq 4$. 
        
    Alice's strategy is as follows. She starts by colouring $r$. In the second round, she colours an uncoloured vertex in $N(u)$ in $G'$ if possible, and if not, then $uv\in E(G)$ and Bob coloured $N(u)\setminus \{u\}$ in the first round, and so, she colours $v$. In the third round, Alice colours an uncoloured vertex in $N(v)$ in $G'$ (this is possible since $d(v)\geq 4$). Then, by the end of the third round, Alice has coloured $r$, at least one vertex in $N(v)$ in $G'$, and at least one vertex in $N(u)$, and these vertices form a connected dominating set of $G$.

    \item $d(u)=1$ and $d(v)\geq 5$.
    
    Let $u'\in N(u)$ be a fixed neighbour of $u$ in $N(u)$. In this case, there exists at least one vertex $r'\in G'$ with $d_{G'}(r')\geq n-5$, as otherwise, we have a contradiction since $$m\leq \frac{(n-3)(n-6)}{2}+2(n-3)+1=\frac{n^2-5n+8}{2}<x.$$ 
    Note that $v$ has at least $4$ neighbours in $G'$ since $d(v)\geq 5$ and $rv\notin E(G)$. We distinguish the following subcases:
    
    \begin{enumerate}
        \item $\Delta(G')=n-4$.
        
        Let $r'\in V(G')$ be such that $d_{G'}(r')=n-4$, then Alice's strategy is as follows. She starts by colouring $u'$ (it may be that $u'=v$). If Bob colours a vertex in $\{r,r'\}$ (a neighbour $v'\in V(G')$ of $v$, resp.) in the first round, then, in the second round, Alice colours the other vertex in $\{r,r'\}$ (another neighbour $v^*\in V(G')$ of $v$, resp.). If Alice coloured a vertex in $\{r,r'\}$ ($v^*$, resp.) in the second round, then she colours a vertex in $\{v',v^*\}$ ($\{r,r'\}$, resp.) in the third round. After three rounds, Alice's vertices form a connected dominating set of $G$.
        
        \item $\Delta(G')=n-5$.
        
        Let $r'\in V(G')$ be such that $d_{G'}(r')=n-5$, and let $G''=G'[N(r')]$. We distinguish cases according to whether $u'\in V(G'')$ or not.
    
            \begin{enumerate}
			\item $u'\in V(G'')$. 
			
			As $u'\in V(G'')$, $uv\notin E(G)$. Since $d(r')=n-5$, there is exactly one additional vertex $w\in V(G')\setminus V(G'')$ ($w\neq r'$). Note that $d_{G'}(w)\geq 1$ because if $d_{G'}(w)=0$, {\it i.e.}, $N(w)=r$, then we have a contradiction since $$m\leq \frac{(n-4)(n-5)}{2}+2(n-3)+1=\frac{n^2-5n+10}{2}=x.$$
			If $d_{G'}(w)\geq 2$, then Alice's strategy is as follows. 
			She starts by colouring $u'$. As before, Bob is forced to colour $r$ in the first round. Indeed, if he does not, then Alice will colour $r$ in the second round, and then she will colour an uncoloured neighbour $v'\in V(G')$ of $v$ in the third round (this is possible since $d(v)\geq 5$), and her vertices form a connected dominating set of $G$. In the second round, Alice colours $r'$. In the third round, if needed, {\it i.e.}, if $u'\notin N(w)$, Alice colours an uncoloured neighbour $w'\in V(G')$ of $w$ (this is possible since $d_{G'}(w)\geq 2$). In the fourth round, Alice colours an uncoloured neighbour $v'\in V(G')$ of $v$ (this is possible since $v$ has at least $5$ neighbours in $G'$ as $uv\notin E(G)$ and $d(v)\geq 5$). 
			At the end of the fourth round, Alice's vertices form a connected dominating set of $G$. 
        
			Assume now that $d_{G'}(w)=1$, and let $N_{G'}(w)=\{w'\}$. Then, $\Delta(G'')=n-6$. Indeed, if $\Delta(G'')\leq n-7$ (and so, $n\geq 7$), then we have a contradiction since $$m\leq \frac{(n-5)(n-7)}{2}+n-5+2(n-3)+1+1=\frac{n^2-6n+17}{2}\leq x.$$ 
			Let $r''\in V(G'')$ be such that $d_{G''}(r'')=n-6$, and observe that $w'\in N(r'')$. Alice's strategy is as follows. She starts by colouring $u'$. As before (when $d_{G'}(w)\geq 2$), Bob is forced to colour $r$ in the first round. In the second round, Alice colours $w'$. Analogously to why Bob was forced to colour $r$ in the first round, Bob is forced to colour $r'$ in the second round. In the third round, Alice colours $r''$. In the fourth round, Alice colours an uncoloured neighbour $v'\in V(G')$ of $v$ (this is possible since $v$ has at least $5$ neighbours in $G'$). 
			At the end of the fourth round, Alice's vertices form a connected dominating set of $G$.   
        
			\item $u'\notin V(G'')$. 
			
			First, assume that $uv\notin E(G)$. Then, $u'$ is the only vertex of $G'$ that is not a neighbour of $r'$, and $d(u')\geq 3$. Indeed, if $d(u')\leq 2$, then we have a contradiction since $$m\leq \frac{(n-4)(n-5)}{2}+2(n-3)+1=\frac{n^2-5n+10}{2}=x.$$
			Thus, there is at least one edge $u'u''$ with $u''\in V(G'')\cup \{v\}$. If $d(u')\geq 4$, then Alice's strategy is as follows. She starts by colouring $u'$, forcing Bob to colour $r$, as before. Then, she colours $r'$, and in the third round, she colours one of the remaining uncoloured neighbours of $u'$ in $G''$ (which exists since $d(u')\geq 4$). In the fourth round, Alice colours an uncoloured neighbour $v'\in V(G')$ of $v$ (this is possible since $v$ has at least $5$ neighbours in $G'$). At the end of the fourth round, Alice's vertices form a connected dominating set of $G$. 
			
			Otherwise, $d(u')=3$, and, as in Case 2.(c)iiA, there exists $r''\in V(G'')$ such that $d_{G''}(r'')=n-6$. Let $r''\in V(G'')$ be such that $d_{G''}(r'')=n-6$. Alice's strategy is as follows. She starts by colouring $u'$, forcing Bob to colour $r$, as before. Then, she colours $u''$, forcing Bob to colour $r'$, as before. In the third round, Alice colours $r''$. In the fourth round, Alice colours an uncoloured neighbour $v'\in V(G')$ of $v$ (this is possible since $v$ has at least $5$ neighbours in $G'$). At the end of the fourth round, Alice's vertices form a connected dominating set of $G$.
			
			Now, assume that $uv\in E(G)$. Then, $u'=v$ and there is exactly one additional vertex $w\in V(G')\setminus V(G'')$ ($w\neq r'$). Note that $d_{G'}(w)\geq 2$ because if $d_{G'}(w)=1$, then we have a contradiction since $$m\leq \frac{(n-4)(n-5)}{2}+2(n-3)+1=\frac{n^2-5n+10}{2}=x.$$
			Alice's strategy is as follows. She starts by colouring $u'=v$, forcing Bob to colour $r$, as before. In the second round, she colours $r'$. In the third round, Alice colours an uncoloured neighbour $w'\in V(G')$ of $w$ (this is possible since $d_{G'}(w)\geq 2$). In the fourth round, if needed, {\it i.e.}, if Alice has not yet coloured a vertex in $N(v)$ that is not $u$, Alice colours an uncoloured neighbour $v'\in V(G')$ of $v$ (this is possible since $d(v)\geq 5$). At the end of the fourth round, Alice's vertices form a connected dominating set of $G$. \qedhere
		    \end{enumerate}
    \end{enumerate}
    \end{enumerate}
\end{enumerate}
\end{proof}

We note that the bound in the statement of Theorem~\ref{thm:dense} is sharp, in the sense that there exists a graph $G$ with  $\frac{(|V(G)|-2)(|V(G)|-3)}{2}+2$ edges that is not A-perfect. For example, consider, as $G$, any graph obtained from a complete graph on an odd number $N \geq 3$ of vertices,
by taking any of its vertices $u$, and attaching at $u$ a pending path $(u,v,w)$ of length~$2$. Note that $|V(G)|=N+2$ and that $$|E(G)|=\frac{N(N-1)}{2}+2=\frac{(|V(G)|-2)(|V(G)|-3)}{2}+2.$$
Now, to see that $G$ is not A-perfect, consider the following strategy for Bob. Bob colours a vertex in $\{u,v\}$ in the first round,
and then, in each of the subsequent rounds, he colours any uncoloured vertex different from $w$. Since $|V(G)|$ is odd, Alice is forced to colour $w$ at some point, which, by the end of the game, cannot be part of a single connected red component due to Bob having coloured $u$ or $v$ in the first round. Thus, $G$ is not A-perfect.

Since the vertices' degrees influence a graph being A-perfect or not, we study regular graphs next, as they are a special case since all of their vertices have the same degree. We prove that there exist arbitrarily large connected $d$-regular graphs $G$, with $d \geq 3$, for which $c_g(G)$ is close to the lower bound (Lemma~\ref{theorem regular d}), while, for every $d \geq 4$, there exist arbitrarily large connected $d$-regular graphs $G$ that are A-perfect (Lemma~\ref{theorem regular n}). However, the latter result does not hold for every $d \geq 3$ since we prove that any sufficiently large cubic graph is not A-perfect (Theorem~\ref{theorem:cubic-main-result}).

Before starting, let us mention the case of $2$-regular graphs, {\it i.e.}, cycles. From a result for paths in~\cite{paper1}, it follows that, for every $n \geq 3$, $c_g(C_n) = 2$ (the lower bound is trivial, and for the upper bound, Bob's strategy is to colour a vertex adjacent to the red vertex in the first round, and now, the game is equivalent to one on a path $P_{n-1}$, with one of its ends initially coloured red). We now show that the lower bound in Lemma~\ref{lemma bounds} is almost tight for arbitrarily large connected $d$-regular graphs, for every $d \geq 3$.

\begin{lemma} \label{theorem regular d}
For every $d\geq 3$, there exist arbitrarily large connected $d$-regular graphs $G$ such that $c_g(G) \leq \left\lceil\frac{d+3}{2}\right\rceil$.
\end{lemma}

\begin{proof}
Let $G$ be the graph constructed as follows. 
Start from $N\geq 2$ disjoint copies $H_0,\dots,H_{N-1}$ of the complete graph on $d+1$ vertices. Now, for every $i \in \{0, \dots, N-1\}$, remove the edge $u_iv_i$, where $u_i$ and $v_i$ are any two distinct vertices of $H_i$. Finally, add the edge $v_iu_{i+1}$ for every $i \in \{0, \dots, N-1\}$ (where, here and further, operations are understood modulo~$N$).
Note that the resulting graph $G$ is $d$-regular, and, free to consider large values of $N$, can be as large as desired.
For every $i \in \{0, \dots, N-1\}$, every vertex of $H_i$ that is different from $u_i$ and $v_i$ is said to be {\it internal} (to $H_i$). Since $d \geq 3$, every $H_i$ has at least two internal vertices.

We give a strategy for Bob that ensures that Alice's score in $G$ is at most $\left\lceil\frac{d+3}{2}\right\rceil$. In each round, if the last vertex coloured by Alice is 

\begin{itemize}
	\item some vertex $u_i$, then Bob colours $v_{i-1}$; 
	\item some vertex $v_i$, then Bob colours $u_{i+1}$;
 	\item a vertex internal to some $H_i$, then Bob colours an uncoloured vertex internal to the same $H_i$.
\end{itemize}
	
	By this strategy, once the game ends, every connected red component must be completely contained inside some $H_i$. This is because this strategy guarantees that any two vertices $v_i$ and $u_{i+1}$ end up coloured either by different players, or by Bob only.
	It thus follows that the largest connected red component contains, in the worst-case scenario, some $u_i$, $v_i$, 
	and half of the other vertices of $H_i$. 
	In other words, the largest connected red component is of order at most $\left\lceil \frac{d+3}{2}\right\rceil$.
\end{proof}

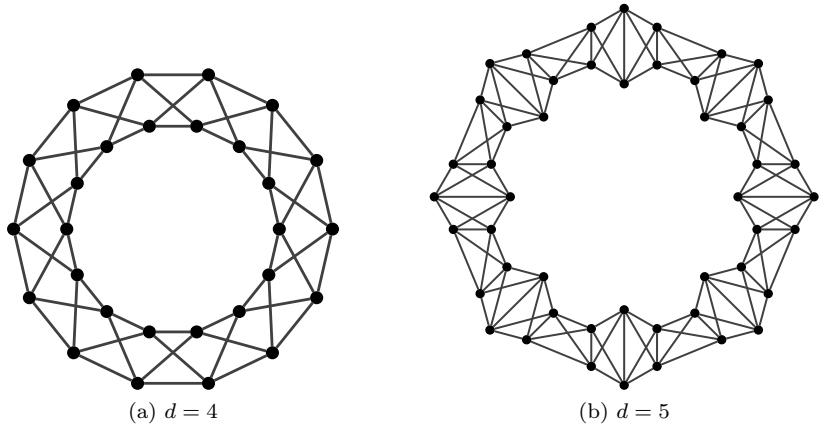
\begin{figure}
    \centering
    \subfloat[$d=4$]{
    \scalebox{0.7}{
	\begin{tikzpicture}
		\draw ( 3.0 , 0.0 ) node[v](l0c1){} ;
		\draw ( 2.7029070971425466 , 1.3016501927232413 ) node[v](l1c1){} ;
		\draw ( 1.870471183855699 , 2.3454930292724088 ) node[v](l2c1){} ;
		\draw ( 0.6675661280869136 , 2.9247829773559344 ) node[v](l3c1){} ;
		\draw ( -0.6675583669103069 , 2.9247847487923013 ) node[v](l4c1){} ;
		\draw ( -1.8704649598727507 , 2.345497992727393 ) node[v](l5c1){} ;
		\draw ( -2.7029036430873643 , 1.3016573651253445 ) node[v](l6c1){} ;
		\draw ( -2.999999999989438 , 7.960769380058192e-06 ) node[v](l7c1){} ;
		\draw ( -2.702910551178696 , -1.3016430203119738 ) node[v](l8c1){} ;
		\draw ( -1.8704774078254758 , -2.345488065800909 ) node[v](l9c1){} ;
		\draw ( -0.6675738892588192 , -2.9247812058989737 ) node[v](l10c1){} ;
		\draw ( 0.6675506057290013 , -2.9247865202080723 ) node[v](l11c1){} ;
		\draw ( 1.8704587358766303 , -2.3455029561658622 ) node[v](l12c1){} ;
		\draw ( 2.7029001890131488 , -1.3016645375182825 ) node[v](l13c1){} ;
		
		\draw ( 2.0 , 0.0 ) node[v](l0c2){} ;
		\draw ( 1.8019380647616978 , 0.8677667951488276 ) node[v](l1c2){} ;
		\draw ( 1.2469807892371327 , 1.563662019514939 ) node[v](l2c2){} ;
		\draw ( 0.44504408539127577 , 1.9498553182372897 ) node[v](l3c2){} ;
		\draw ( -0.4450389112735379 , 1.9498564991948675 ) node[v](l4c2){} ;
		\draw ( -1.2469766399151672 , 1.5636653284849287 ) node[v](l5c2){} ;
		\draw ( -1.8019357620582428 , 0.8677715767502296 ) node[v](l6c2){} ;
		\draw ( -1.9999999999929585 , 5.307179586705461e-06 ) node[v](l7c2){} ;
		\draw ( -1.801940367452464 , -0.8677620135413159 ) node[v](l8c2){} ;
		\draw ( -1.2469849385503171 , -1.5636587105339395 ) node[v](l9c2){} ;
		\draw ( -0.4450492595058795 , -1.9498541372659823 ) node[v](l10c2){} ;
		\draw ( 0.44503373715266753 , -1.9498576801387149 ) node[v](l11c2){} ;
		\draw ( 1.2469724905844202 , -1.5636686374439082 ) node[v](l12c2){} ;
		\draw ( 1.8019334593420993 , -0.8677763583455217 ) node[v](l13c2){} ;
		
		\Edge[](l1c1)(l0c2)
		\Edge[](l1c1)(l0c1)
		\Edge[](l1c2)(l0c2)
		\Edge[](l1c2)(l0c1)
		\Edge[](l13c1)(l0c2)
		\Edge[](l13c1)(l0c1)
		\Edge[](l13c2)(l0c2)
		\Edge[](l13c2)(l0c1)
		\Edge[](l1c1)(l2c2)
		\Edge[](l1c1)(l2c1)
		\Edge[](l1c2)(l2c2)
		\Edge[](l1c2)(l2c1)
		\Edge[](l2c1)(l3c2)
		\Edge[](l2c1)(l3c1)
		\Edge[](l2c2)(l3c2)
		\Edge[](l2c2)(l3c1)
		\Edge[](l3c1)(l4c2)
		\Edge[](l3c1)(l4c1)
		\Edge[](l3c2)(l4c2)
		\Edge[](l3c2)(l4c1)
		\Edge[](l4c1)(l5c2)
		\Edge[](l4c1)(l5c1)
		\Edge[](l4c2)(l5c2)
		\Edge[](l4c2)(l5c1)
		\Edge[](l5c1)(l6c2)
		\Edge[](l5c1)(l6c1)
		\Edge[](l5c2)(l6c2)
		\Edge[](l5c2)(l6c1)
		\Edge[](l6c1)(l7c2)
		\Edge[](l6c1)(l7c1)
		\Edge[](l6c2)(l7c2)
		\Edge[](l6c2)(l7c1)
		\Edge[](l7c1)(l8c2)
		\Edge[](l7c1)(l8c1)
		\Edge[](l7c2)(l8c2)
		\Edge[](l7c2)(l8c1)
		\Edge[](l8c1)(l9c2)
		\Edge[](l8c1)(l9c1)
		\Edge[](l8c2)(l9c2)
		\Edge[](l8c2)(l9c1)
		\Edge[](l9c1)(l10c2)
		\Edge[](l9c1)(l10c1)
		\Edge[](l9c2)(l10c2)
		\Edge[](l9c2)(l10c1)
		\Edge[](l10c1)(l11c2)
		\Edge[](l10c1)(l11c1)
		\Edge[](l10c2)(l11c2)
		\Edge[](l10c2)(l11c1)
		\Edge[](l11c1)(l12c2)
		\Edge[](l11c1)(l12c1)
		\Edge[](l11c2)(l12c2)
		\Edge[](l11c2)(l12c1)
		\Edge[](l12c1)(l13c2)
		\Edge[](l12c1)(l13c1)
		\Edge[](l12c2)(l13c2)
		\Edge[](l12c2)(l13c1)
	\end{tikzpicture}
    }
    }
    \hspace{20pt}
    \subfloat[$d=5$]{
    \scalebox{0.5}{
    \begin{tikzpicture}
		\draw ( 5.0 , 0.0 ) node[v](h0s0){} ;
		\draw ( 4.500000766025195 , 0.8660249615191342 ) node[v](h0s1){} ;
		\draw ( 3.500001532051564 , 0.8660262883130148 ) node[v](h0s2){} ;
		\draw ( 3.0000000000035207 , 2.6535897937968196e-06 ) node[v](h0s3){} ;
		\draw ( 3.4999969359015664 , -0.8660236347191557 ) node[v](h0s4){} ;
		\draw ( 4.499996169868156 , -0.8660276151007966 ) node[v](h0s5){} ;
		\Edge[](h0s1)(h0s0)
		\Edge[](h0s2)(h0s0)
		\Edge[](h0s2)(h0s1)
		\Edge[](h0s3)(h0s0)
		\Edge[](h0s3)(h0s1)
		\Edge[](h0s3)(h0s2)
		\Edge[](h0s4)(h0s0)
		\Edge[](h0s4)(h0s3)
		\Edge[](h0s5)(h0s0)
		\Edge[](h0s5)(h0s3)
		\Edge[](h0s5)(h0s4)
		\draw ( 3.5355362513961315 , 3.535531560467788 ) node[v](h1s0){} ;
		\draw ( 2.5696114511977277 , 3.7943514752952274 ) node[v](h1s1){} ;
		\draw ( 1.8625038043969726 , 3.087246643049323 ) node[v](h1s2){} ;
		\draw ( 2.121319874470075 , 2.121320812655745 ) node[v](h1s3){} ;
		\draw ( 3.0872439878631943 , 1.8624983346610247 ) node[v](h1s4){} ;
		\draw ( 3.794353511027625 , 2.569601290530806 ) node[v](h1s5){} ;
		\Edge[](h1s1)(h1s0)
		\Edge[](h1s2)(h1s0)
		\Edge[](h1s2)(h1s1)
		\Edge[](h1s3)(h1s0)
		\Edge[](h1s3)(h1s1)
		\Edge[](h1s3)(h1s2)
		\Edge[](h1s4)(h1s0)
		\Edge[](h1s4)(h1s3)
		\Edge[](h1s5)(h1s0)
		\Edge[](h1s5)(h1s3)
		\Edge[](h1s5)(h1s4)
		\draw ( 6.633974483387665e-06 , 4.999999999995599 ) node[v](h2s0){} ;
		\draw ( -0.8660189909403205 , 4.500001915058733 ) node[v](h2s1){} ;
		\draw ( -0.8660216445280814 , 3.5000026810877434 ) node[v](h2s2){} ;
		\draw ( 1.3267948962427866e-06 , 3.000000000004401 ) node[v](h2s3){} ;
		\draw ( 0.8660282784964664 , 3.4999957868627467 ) node[v](h2s4){} ;
		\draw ( 0.8660335856719876 , 4.499995020823176 ) node[v](h2s5){} ;
		\Edge[](h2s1)(h2s0)
		\Edge[](h2s2)(h2s0)
		\Edge[](h2s2)(h2s1)
		\Edge[](h2s3)(h2s0)
		\Edge[](h2s3)(h2s1)
		\Edge[](h2s3)(h2s2)
		\Edge[](h2s4)(h2s0)
		\Edge[](h2s4)(h2s3)
		\Edge[](h2s5)(h2s0)
		\Edge[](h2s5)(h2s3)
		\Edge[](h2s5)(h2s4)
		\draw ( -3.5355268695332196 , 3.5355409423182516 ) node[v](h3s0){} ;
		\draw ( -3.7943480659445274 , 2.5696164855216406 ) node[v](h3s1){} ;
		\draw ( -3.087244171886063 , 1.8625079005384246 ) node[v](h3s2){} ;
		\draw ( -2.121317998097494 , 2.121322689025837 ) node[v](h3s3){} ;
		\draw ( -1.862494238519817 , 3.0872464590137625 ) node[v](h3s4){} ;
		\draw ( -2.5695962561996692 , 3.794356920358165 ) node[v](h3s5){} ;
		\Edge[](h3s1)(h3s0)
		\Edge[](h3s2)(h3s0)
		\Edge[](h3s2)(h3s1)
		\Edge[](h3s3)(h3s0)
		\Edge[](h3s3)(h3s1)
		\Edge[](h3s3)(h3s2)
		\Edge[](h3s4)(h3s0)
		\Edge[](h3s4)(h3s3)
		\Edge[](h3s5)(h3s0)
		\Edge[](h3s5)(h3s3)
		\Edge[](h3s5)(h3s4)
		\draw ( -4.999999999982396 , 1.3267948966763651e-05 ) node[v](h4s0){} ;
		\draw ( -4.500003064084351 , -0.8660130203599824 ) node[v](h4s1){} ;
		\draw ( -3.5000038301177607 , -0.8660170007416236 ) node[v](h4s2){} ;
		\draw ( -3.0 , 5.307179586280058e-06 ) node[v](h4s3){} ;
		\draw ( -3.4999946378177658 , 0.8660329222722525 ) node[v](h4s4){} ;
		\draw ( -4.499993871770273 , 0.8660395562416541 ) node[v](h4s5){} ;
		\Edge[](h4s1)(h4s0)
		\Edge[](h4s2)(h4s0)
		\Edge[](h4s2)(h4s1)
		\Edge[](h4s3)(h4s0)
		\Edge[](h4s3)(h4s1)
		\Edge[](h4s3)(h4s2)
		\Edge[](h4s4)(h4s0)
		\Edge[](h4s4)(h4s3)
		\Edge[](h4s5)(h4s0)
		\Edge[](h4s5)(h4s3)
		\Edge[](h4s5)(h4s4)
		\draw ( -3.5355456332341455 , -3.53552217859243 ) node[v](h5s0){} ;
		\draw ( -2.569621519841027 , -3.7943446565871497 ) node[v](h5s1){} ;
		\draw ( -1.862511996676596 , -3.0872417007173687 ) node[v](h5s2){} ;
		\draw ( -2.121325503577863 , -2.1213151835355095 ) node[v](h5s3){} ;
		\draw ( -3.0872489301588955 , -1.8624901423753324 ) node[v](h5s4){} ;
		\draw ( -3.794360329682023 , -2.5695912218640107 ) node[v](h5s5){} ;
		\Edge[](h5s1)(h5s0)
		\Edge[](h5s2)(h5s0)
		\Edge[](h5s2)(h5s1)
		\Edge[](h5s3)(h5s0)
		\Edge[](h5s3)(h5s1)
		\Edge[](h5s3)(h5s2)
		\Edge[](h5s4)(h5s0)
		\Edge[](h5s4)(h5s3)
		\Edge[](h5s5)(h5s0)
		\Edge[](h5s5)(h5s3)
		\Edge[](h5s5)(h5s4)
		\draw ( -1.9901923452336725e-05 , -4.999999999960391 ) node[v](h6s0){} ;
		\draw ( 0.8660070497781176 , -4.500004213102046 ) node[v](h6s1){} ;
		\draw ( 0.8660123569536394 , -3.5000049791416177 ) node[v](h6s2){} ;
		\draw ( -9.28756427764025e-06 , -2.9999999999903175 ) node[v](h6s3){} ;
		\draw ( -0.8660375660465156 , -3.499993488766623 ) node[v](h6s4){} ;
		\draw ( -0.8660455268097978 , -4.499992722709448 ) node[v](h6s5){} ;
		\Edge[](h6s1)(h6s0)
		\Edge[](h6s2)(h6s0)
		\Edge[](h6s2)(h6s1)
		\Edge[](h6s3)(h6s0)
		\Edge[](h6s3)(h6s1)
		\Edge[](h6s3)(h6s2)
		\Edge[](h6s4)(h6s0)
		\Edge[](h6s4)(h6s3)
		\Edge[](h6s5)(h6s0)
		\Edge[](h6s5)(h6s3)
		\Edge[](h6s5)(h6s4)
		\draw ( 3.5355174876454134 , -3.535550324143819 ) node[v](h7s0){} ;
		\draw ( 3.7943412472230893 , -2.5696265541558936 ) node[v](h7s1){} ;
		\draw ( 3.087239229543238 , -1.8625160928114917 ) node[v](h7s2){} ;
		\draw ( 2.1213123689697886 , -2.1213283181261566 ) node[v](h7s3){} ;
		\draw ( 1.862486046227566 , -3.0872514012985945 ) node[v](h7s4){} ;
		\draw ( 2.569586187523825 , -3.7943637389992046 ) node[v](h7s5){} ;
		\Edge[](h7s1)(h7s0)
		\Edge[](h7s2)(h7s0)
		\Edge[](h7s2)(h7s1)
		\Edge[](h7s3)(h7s0)
		\Edge[](h7s3)(h7s1)
		\Edge[](h7s3)(h7s2)
		\Edge[](h7s4)(h7s0)
		\Edge[](h7s4)(h7s3)
		\Edge[](h7s5)(h7s0)
		\Edge[](h7s5)(h7s3)
		\Edge[](h7s5)(h7s4)
		\Edge[](h0s1)(h1s4)
		\Edge[](h0s2)(h1s4)
		\Edge[](h0s1)(h1s5)
		\Edge[](h0s2)(h1s5)
		\Edge[](h1s1)(h2s4)
		\Edge[](h1s2)(h2s4)
		\Edge[](h1s1)(h2s5)
		\Edge[](h1s2)(h2s5)
		\Edge[](h2s1)(h3s4)
		\Edge[](h2s2)(h3s4)
		\Edge[](h2s1)(h3s5)
		\Edge[](h2s2)(h3s5)
		\Edge[](h3s1)(h4s4)
		\Edge[](h3s2)(h4s4)
		\Edge[](h3s1)(h4s5)
		\Edge[](h3s2)(h4s5)
		\Edge[](h4s1)(h5s4)
		\Edge[](h4s2)(h5s4)
		\Edge[](h4s1)(h5s5)
		\Edge[](h4s2)(h5s5)
		\Edge[](h5s1)(h6s4)
		\Edge[](h5s2)(h6s4)
		\Edge[](h5s1)(h6s5)
		\Edge[](h5s2)(h6s5)
		\Edge[](h6s1)(h7s4)
		\Edge[](h6s2)(h7s4)
		\Edge[](h6s1)(h7s5)
		\Edge[](h6s2)(h7s5)
		\Edge[](h7s1)(h0s4)
		\Edge[](h7s2)(h0s4)
		\Edge[](h7s1)(h0s5)
		\Edge[](h7s2)(h0s5)
	\end{tikzpicture}
    }
    }
    \caption{Examples of $d$-regular A-perfect graphs constructed in the proof of Lemma~\ref{theorem regular n}.}
    \label{figure:aperfect-regular}
\end{figure}

Regarding the upper bound from Lemma~\ref{lemma bounds} in the context of arbitrarily large connected $d$-regular graphs, we prove the following:

\begin{lemma} \label{theorem regular n}
For every $d \ge 4$, there exist arbitrarily large  $d$-regular A-perfect graphs.
\end{lemma}

\begin{proof}
Let $N > 2$, and let $d \geq 4$ be fixed.
To prove the claim, we construct a $d$-regular graph $G$, whose order is a function of $N$,
such that $c_g(G)=\left\lceil \frac{|V(G)|}{2} \right\rceil$.
We give two possible constructions for $G$, depending on whether $d=4$ or $d \geq 5$ (see Figure~\ref{figure:aperfect-regular} for an illustration of both constructions).

\begin{itemize}
	\item For the case $d=4$, $G$ is the $4$-regular graph having two vertices $u^i_1$ and $u^i_2$ for every $i \in \{0, \dots, N-1\}$,
	and the four edges $u^i_1u^{i+1}_1$, $u^i_1u^{i+1}_2, u^i_2u^{i+1}_1$, $u^i_2u^{i+1}_2$ for every $i \in \{0, \dots, N-1\}$ (where, here and further,
	operations over the superscripts are modulo~$N$).
	
	To prove that $G$ is A-perfect, we give a strategy for Alice that ensures that, at the end of the game in $G$, the red subgraph is connected. In the first round, Alice colours $u_1^0$. Then, in the subsequent rounds, if the last vertex Bob coloured is $u_1^j$ ($u_2^j$, resp.) for some $j \in \{1, \dots, N-1\}$, Alice responds by colouring $u_2^j$ ($u_1^j$, resp.). Otherwise, Alice colours an arbitrary uncoloured vertex. By Alice's strategy, at the end of the game, for every $0\leq i\leq N-1$, exactly one of $u_1^i$ and $u_2^i$ is red, and thus, the red subgraph is connected, and $G$ is A-perfect. 
	
	\item We now consider the case where $d \geq 5$. Here, $G$ is constructed as follows.
	Start from $N$ disjoint copies $H_0, \dots, H_{N-1}$ of the complete graph on $d+1$ vertices, where, for every $i \in \{0,\dots,N-1\}$, we denote by $v_1^i,\dots,v_{d+1}^i$ the vertices of $H_i$.
	For every $i \in \{0,\dots,N-1\}$, we remove the edges $v^i_1v^i_3$, $v^i_1v^i_4$, $v^i_2v^i_3$ and $v^i_2v^i_4$ from $H_i$.
	To finish the construction of $G$ and make it $d$-regular, we then join the $H_i$'s by adding the edges $v^i_3v^{i+1}_1$, $v^i_3v^{i+1}_2$, $v^i_4v^{i+1}_1$, and $v^i_4v^{i+1}_2$ for every $i \in \{0,\dots,N-1\}$ 
	(again, operations are understood modulo $N$).
	
	To prove that $G$ is A-perfect, we give a strategy for Alice that ensures her a score of $\left \lceil |V(G)|/2 \right \rceil$. In the first round, Alice colours any vertex. In each of the subsequent rounds, if the last vertex Bob coloured is
	
	\begin{itemize}
	    \item in some pair $\{v^i_1, v^i_2\}$ or $\{v^i_3, v^i_4\}$, then Alice colours the other vertex in that pair;
	    \item some vertex $v^i_j$ with $5\leq j\leq d+1$, then Alice colours another vertex $v^i_{\ell}$ with $5\leq \ell \leq d+1$ and $j\neq \ell$. 
	\end{itemize}
	
	Whenever Alice cannot follow the strategy above, she colours an arbitrary uncoloured vertex. By Alice's strategy, at the end of the game, for every $i \in \{0,\dots,N-1\}$, at least one vertex in $\{v^i_1, v^i_2\}$ is red, at least one vertex in $\{v^i_3, v^i_4\}$ is red, and at least one vertex in $\{v^i_5,\ldots,v^i_{d+1}\}$ is red. These vertices form a connected dominating set of $G$, from which we deduce that $c_g(G) = \left \lceil |V(G)|/2 \right \rceil$, by Lemma~\ref{lemma connected dominating set}. Thus, $G$ is A-perfect. \qedhere
\end{itemize}
\end{proof}

As mentioned earlier, the bound on $d$ in the statement of Lemma~\ref{theorem regular n} cannot be lowered, as we prove that A-perfect cubic graphs have bounded order. This can actually be established through previous results on the existence of particular cuts in sufficiently large connected cubic graphs, such as ones from~\cite{cubic} relying on the following terminology.

A \textit{supercycle} is a connected graph with minimum degree at least~$2$ where not all vertices are of degree~$2$. For a graph $G$, a matching $M$ is said \textit{suitable} if $G-M$ consists of exactly two connected components, each of which is a supercycle. Note that if $G$ is cubic, then, in $G-M$, every vertex incident to an edge of $M$ has degree precisely~$2$, while, by the definition of a supercycle, each of the two connected components contains a degree-$3$ vertex.

In~\cite{cubic}, Mukkamala and P{\'a}lv{\"o}lgyi proved the following result on the existence of suitable matchings in sufficiently large connected cubic graphs.

\begin{theorem}[Mukkamala and P{\'a}lv{\"o}lgyi, Corollary~1 of~\cite{cubic}]\label{theorem:cubic-suitable-matching}
Every connected cubic graph with order at least~$18$ admits a suitable matching.
\end{theorem}

We are now ready to prove the aforementioned result on cubic graphs.

\begin{theorem}\label{theorem:cubic-main-result}
Every A-perfect cubic graph has order at most $16$.
\end{theorem}

\begin{proof}
Let $G$ be an A-perfect cubic graph. Since $G$ is cubic, each of its connected components has order at least $4$. Thus, by Corollary~\ref{cor:A-perfect}, $G$ is connected. 
Towards proving the claim, assume that $G$ has order at least~$18$.
Then, by Theorem~\ref{theorem:cubic-suitable-matching}, $G$ admits a suitable matching $M$. As mentioned earlier, $G-M$ consists of exactly two connected components $C_1$ and $C_2$, in each of which all vertices incident to an edge of $M$ have degree exactly~$2$ while the other vertices (there is at least one such) have degree exactly~$3$. Since, by the handshaking lemma, in every graph the number of odd-degree vertices is even, we deduce that, in each of $C_1$ and $C_2$, there are actually at least two degree-$3$ vertices. In what follows, the degree-$2$ vertices of the $C_i$'s are called \textit{interior vertices}, while their other (degree-$3$) vertices are called \textit{exterior vertices}.

Consider now the strategy for Bob where, each turn during a game on $G$, he answers to Alice's moves as follows:

\begin{itemize}
    \item if Alice colours an interior vertex incident to an edge $e \in M$, then Bob colours the second interior vertex incident to $e$;
    
    \item if Alice colours an exterior vertex $v$, then Bob plays as follows:
    
    \begin{itemize}
        \item if $v$ is the first exterior vertex coloured by Alice during the whole game, then, denoting, for the rest of the game, by $C^*$ the one of $C_1$ and $C_2$ that contains $v$, Bob colours any uncoloured exterior vertex of $C^*$ (one such exists, since $C^*$ contains at least two exterior vertices);
        
        \item if $v$ is not the first exterior vertex that Alice colours, then $C^*$ was defined during an earlier turn, and Bob colours any uncoloured exterior vertex of $C^*$. If $C^*$ does not contain any such uncoloured vertex, then Bob colours any uncoloured vertex of $G$ instead.
    \end{itemize}
\end{itemize}

Note that Bob can clearly follow the strategy above from start to end.
Once the game ends, note also that for every edge of $M$ the two incident vertices are coloured with distinct colours. Furthermore, since $C_1$ and $C_2$ have at least two exterior vertices each, each $C_i$ must contain an exterior vertex $u_i$ coloured red. From all these arguments, we deduce that the red subgraph cannot contain a path joining $u_1$ and $u_2$, and thus the red subgraph is not connected. Thus, an A-perfect cubic graph must have order strictly less than~$18$.
\end{proof}

\section{Graphs with few $P_4$'s}\label{section:sparse}

In this section, we give a linear-time algorithm determining $c_g$ for $(q,q-4)$-graphs, which are graphs containing few $P_4$'s~\cite{P4sparse}. For a fixed $q\geq 0$, a graph $G$ is a {\it $(q,q-4)$-graph} if every subset $S \subseteq V(G)$ of at most $q$ vertices induces at most $q-4$ paths on $4$ vertices. Note that a cograph is a $(q,q-4)$-graph when $q=4$.

A better understanding of our game in these graphs would allow to establish a similar result for the Largest Connected Subgraph game, for which a linear-time algorithm is only known for cographs~\cite{paper1}. Recall that, for two graphs $G$ and $H$, $G + H$ is the {\it disjoint union} of $G$ and $H$, where $V(G + H) = V(G) \cup V(H)$ and $E(G + H) = E(G) \cup E(H)$, and $G \oplus H$ is the {\it join} of $G$ and $H$, where $V(G \oplus H) = V(G) \cup V(H)$ and $E(G \oplus H) = E(G) \cup E(H) \cup  \left \{uv : (u,v) \in V(G)\times V(H)\right \}$. Let us give a characterisation of $(q,q-4)$ graphs. Let $G=(S,K,R,E)$ be a graph with $V(G)=S\cup K\cup R$ and $E(G)=E$. Consider the following properties:
\begin{enumerate}
 \item $S \cup K \cup R$ is a partition of $V(G)$, where $R$ can be empty. 
 \item $G[K\cup R]= K\oplus R$ ({\it i.e.}, $uv\in E$ for all $u,v\in V(G)$ with $u\in K$ and $v\in R$), and $K$ separates $S$ from $R$ ({\it i.e.}, $uv\notin E$ for all $u\in S$ and $v\in R$). 
 \item $S$ is an independent set, $K$ is a clique, $|S|=|K|\geq 2$, and there exists a bijection $f: S\rightarrow K$ such that, for all $s\in S$, either $N(s)\cap K=K\setminus \{f(s)\}$ or $N(s)\cap K=\{f(s)\}$. In the former case, we say that $f$ is an {\it antimatching}, with the vertices $s$ and $f(s)$ being {\it antimatched}, and in the latter case, we say that $f$ is a {\it matching}, with the vertices $s$ and $f(s)$ being {\it matched}.
\end{enumerate}

If $G=(S,K,R,E)$ verifies all the properties above, it is called a {\it spider}. In that case, if $f$ is a matching (antimatching, resp.), then $G$ is a {\it matched spider} ({\it antimatched spider}, resp.). Also, if $G$ only verifies Properties 1 and 2 above, it is called a {\it pseudo-spider}. In this case, for any fixed $q\geq 0$ such that $|V(S\cup K)|\leq q$, $G$ is a {\it $q$-pseudo-spider}.  A graph $G$ is a $(q,q-4)$-graph if $G$ is the graph $K_1$ or one of the following is satisfied~\cite{Babel95}: 
\begin{enumerate}
	\item $G= G_1 + G_2$, where $G_1$ and $G_2$ are $(q,q-4)$-graphs. 
		\item  $G = G_1 \oplus G_2$, where $G_1$ and $G_2$ are $(q,q-4)$-graphs.
		\item $G$ is the spider $(S,K,R,E)$, where $G[R]$ (if $R$ is not empty) is a $(q,q-4)$-graph. By the definition of a spider, $G[S \cup K]$ induces a $(q,q-4)$-graph.
		\item $G$ is the $q$-pseudo-spider $(S,K,R,E)$, where $G[R]$ (if $R$ is not empty) is a $(q,q-4)$-graph. 
\end{enumerate}

The above definition is a recursive one, in that, for every $(q,q-4)$-graph $G$, there exists a (not necessarily unique) {\it decomposition-tree} representing $G$. The internal nodes of such a tree correspond to subgraphs of $G$ that are $(q,q-4)$-graphs, and its leaves either correspond to a single vertex or to a subgraph with at most $q$ vertices. The root corresponds to $G$, and every internal node has two children (describing cases $1$ to $4$ above). Such a decomposition-tree can be computed in linear time~\cite{BO99}. We can now prove the main result in this section:

\begin{theorem} \label{theorem:qgraphs}
Let $q \geq 0$. For a $(q,q-4)$-graph $G$, determining $c_g(G)$ and an optimal strategy for Alice can be done in linear time.
\end{theorem}

\begin{proof}
Let us first compute (in linear time) a decomposition-tree $T$ of $G$. Now, let us describe the algorithm that proceeds bottom-up from the leaves to the root of $T$. Every leaf of $T$ corresponds to a subgraph $G'$ with a bounded number of vertices, and therefore, $c_g(G')$ and an optimal strategy for Alice can be computed in time $O(1)$. For every internal node $v$ (corresponding to a subgraph $G'$ of $G$) of $T$, $c_g(G')$ and a corresponding strategy for Alice are computed from what has already  been computed for the two subgraphs corresponding to the children of $v$. Precisely, let $G_1$ and $G_2$ be the two subgraphs of $G$ corresponding to the children of the root of $T$, and assume by induction that $c_g(G_1)$, $c_g(G_2)$, and optimal strategies for Alice in $G_1$ and $G_2$ have been computed in linear time. We now describe how the algorithm proceeds for $G$, and we set $|V(G)|=n$. There are $4$ cases depending on how $G$ is obtained from $G_1$ and $G_2$. 
\begin{enumerate}
	\item $G=G_1+ G_2$. Then, $c_g(G)=\max\{c_g(G_1),c_g(G_2)\}$ by Lemma~\ref{lemma subgraph}. W.l.o.g., $c_g(G)=c_g(G_1)$. By induction, $c_g(G_1)$ and a strategy for Alice have already been computed.
	
	\item $G=G_1\oplus G_2$. Then, $c_g(G)=\left\lceil\frac{n}{2}\right\rceil$.
	Indeed, w.l.o.g., $|V(G_1)|\leq |V(G_2)|$, and Alice obtains a connected dominating set by first colouring a vertex in $G_1$, and then one in $G_2$ (unless $G$ is an edge). The result follows by Lemma~\ref{lemma connected dominating set}.
	
	\item $G=(S,K,R,E)$ is a spider with $G_1=G[S\cup K]$ and $G_2=G[R]$. Note that $|R|$ and $n$ have the same parity since $|S|=|K|$. There are two subcases:
	\begin{enumerate}
		\item $G$ is an antimatched spider. Assume that $|K|\geq 3$ since $G$ is a matched spider if $|K|=2$. Then, $c_g(G)=\left\lceil\frac{n}{2}\right\rceil$. Indeed, Alice's strategy is to colour any two uncoloured vertices $v_1,v_2\in K$ in the first two rounds (this is possible since $|K|\geq 3$). Since $G$ is an antimatched spider, for every vertex $v\in S$, at least one of the edges in $\{vv_1,vv_2\}$ is in $E$. Thus, since $K$ is also a clique and $G[K\cup R]= K\oplus R$, the set $\{v_1,v_2\}$ forms a connected dominating set of $G$, and we get the result by Lemma~\ref{lemma connected dominating set}.
		
		\item $G$ is a matched spider. Let us show that
		\begin{equation*}
        c_g(G) = \begin{cases}
            \left\lceil\frac{n}{2}\right\rceil-\left\lceil\frac{\left\lfloor\frac{|K|}{2}\right\rfloor}{2}\right\rceil &\text{if~}  n \text{~and~} \left \lfloor \frac{|K|}{2}\right \rfloor \text{~are odd.} \\[13pt]
            \left\lceil\frac{n}{2}\right\rceil-\left\lfloor\frac{\left\lfloor\frac{|K|}{2}\right\rfloor}{2}\right\rfloor & \text{otherwise.}
           \end{cases}
        \end{equation*}
		\noindent{\bf Bob's strategy.} We give a strategy for Bob to prove the upper bound on $c_g(G)$ in both cases. Bob first plays exhaustively in $K$ ({\it i.e.}, until every vertex in $K$ is coloured), then, while possible, he colours vertices that are not vertices of $S$ matched to blue vertices in $K$, and finally, he colours the vertices of $S$ matched to blue vertices of $K$. At the end of the game, any red vertex in $S$ that is matched to a blue vertex of $K$ forms a one-vertex connected red component. Let $r_S^*$ be the number of such red vertices. Then, $c_g(G)\leq\left\lceil\frac{n}{2}\right\rceil-r_S^*$. 
		
		We first show that $r_S^*\geq \left\lfloor\frac{\left\lfloor\frac{|K|}{2}\right\rfloor}{2}\right\rfloor$. Let $b_K$ be the number of blue vertices in $K$ once every vertex in $K$ is coloured. Since Bob first exhaustively colours the vertices in $K$, then $b_K\geq \left\lfloor \frac{|K|}{2} \right\rfloor$. Then, while it is possible, Bob colours vertices that are not vertices of $S$ matched to blue vertices in $K$. Consider the first point of the game where no such vertex exists (this can occur after any player's move). Let $r_S \geq 0$ be the number of vertices in $S$ that, at this point, are red and matched to a blue vertex in $K$. Now, Bob colours the uncoloured vertices of $S$ matched to blue vertices, and so, Bob colours at most $\left\lceil \frac{b_K-r_S}{2} \right\rceil$ such vertices. Thus, Alice colours at least $\left\lfloor \frac{b_K-r_S}{2} \right\rfloor$ such vertices. Then, $r_S^* \geq r_S + \left\lfloor \frac{b_K-r_S}{2} \right\rfloor \geq \left\lfloor \frac{b_K}{2} \right\rfloor \geq \left\lfloor\frac{\left\lfloor\frac{|K|}{2}\right\rfloor}{2}\right\rfloor$.
		
		Now, consider the case where $n$ and $\left\lfloor\frac{|K|}{2}\right\rfloor$ are odd. We refine the above analysis to show that, in this case, $r_S^*\geq \left\lceil\frac{\left\lfloor\frac{|K|}{2}\right\rfloor}{2}\right\rceil$. If $b_K> \left\lfloor \frac{|K|}{2}\right \rfloor$, then $\left\lfloor \frac{b_K}{2} \right\rfloor > \left\lfloor\frac{\left\lfloor\frac{|K|}{2}\right\rfloor}{2}\right\rfloor$ (as $\left\lfloor \frac{|K|}{2}\right\rfloor$ is odd), and so, $\left\lfloor \frac{b_K}{2} \right\rfloor \geq \left\lceil\frac{\left\lfloor\frac{|K|}{2}\right\rfloor}{2}\right\rceil$, implying that $r^*_S \geq \left\lceil\frac{\left\lfloor\frac{|K|}{2}\right\rfloor}{2}\right\rceil$. Thus, assume that $b_K = \left\lfloor \frac{|K|}{2} \right\rfloor$, and so, $b_K$ is odd. As $n$ is odd, Alice plays last in $G$. Hence, just before Bob colours his first vertex of $S$ matched to a blue vertex in $K$, there are an even number of such uncoloured vertices remaining. Since $b_K$ is odd, then $r_S \geq 1$. Hence, $r_S^* \geq 1 + \left\lfloor \frac{b_K-1}{2} \right\rfloor = 1 + \left\lfloor \frac{\left\lfloor \frac{|K|}{2} \right\rfloor-1}{2} \right\rfloor \geq \left\lceil\frac{\left\lfloor\frac{|K|}{2}\right\rfloor}{2}\right\rceil$. Thus, we have proved the upper bound on $c_g(G)$ in both cases.
		
		\smallskip
		\noindent{\bf Alice's strategy.} Now, we give a strategy for Alice to prove the lower bound on $c_g(G)$ in both cases. Alice follows the same strategy as Bob above. Let $r_K$ be the number of red vertices in $K$ once all the vertices of $K$ are coloured. Since Alice first exhaustively colours the vertices in $K$, we have that $r_K\geq \left\lceil \frac{|K|}{2} \right\rceil$. 
		Let $b_K=|K|-r_K \leq \left\lfloor \frac{|K|}{2} \right\rfloor$ be the number of blue vertices in $K$ once all the vertices of $K$ are coloured. Let $u_S$ be the number of vertices of $S$ that are matched to blue vertices in $K$. Obviously, $u_S \leq b_K$. Alice's strategy ensures that, at the end of the game, the red vertices induce one connected component $X$ and (if Bob plays optimally) some isolated vertices in $S$ that are matched to blue vertices in $K$. By Alice's strategy, there are at most $\left\lceil \frac{u_S}{2} \right\rceil$ such isolated red vertices. Hence, $|X| \geq \left\lceil \frac{n}{2}\right\rceil - \left\lceil \frac{u_S}{2} \right\rceil  \geq \left\lceil \frac{n}{2}\right\rceil - \left\lceil \frac{b_K}{2} \right\rceil$. Thus, $|X| \geq \left\lceil \frac{n}{2}\right\rceil - \left\lceil \frac{\left\lfloor\frac{|K|}{2}\right\rfloor}{2} \right\rceil$, which matches the upper bound when $n$ and $\left\lfloor\frac{|K|}{2}\right\rfloor$ are odd. Also, if $\left\lfloor\frac{|K|}{2}\right\rfloor$ is even, then $\left\lceil \frac{\left\lfloor\frac{|K|}{2}\right\rfloor}{2} \right\rceil=\left\lfloor \frac{\left\lfloor\frac{|K|}{2}\right\rfloor}{2} \right\rfloor$, and so, $|X|\geq \left\lceil \frac{n}{2}\right\rceil - \left\lfloor \frac{\left\lfloor\frac{|K|}{2}\right\rfloor}{2} \right\rfloor$. 
		
		The last case to consider is when $n$ is even. Then, Bob plays last in $G$. This implies that Alice colours at most $\left\lfloor \frac{u_S}{2} \right\rfloor$ vertices of $S$ matched to blue vertices in $K$. So, $|X| \geq \left\lceil \frac{n}{2}\right\rceil - \left\lfloor \frac{u_S}{2} \right\rfloor  \geq \left\lceil \frac{n}{2}\right\rceil - \left\lfloor \frac{b_K}{2} \right\rfloor \geq \left\lceil \frac{n}{2}\right\rceil - \left\lfloor \frac{\left\lfloor\frac{|K|}{2}\right\rfloor}{2} \right\rfloor.$
	\end{enumerate}
  
	\item Finally, let us assume that $G=(S,K,R,E)$ is a $q$-pseudo-spider with $G_1=G[S\cup K]$ (with $|V(G_1)|\leq q$) and $G_2=G[R]$. By Lemma~\ref{lemma subgraph}, we may assume that $G$ is connected. 
	
	First, let us consider the case where $|V(G_2)|\leq 2q$, and so, $|V(G)|\leq 3q$. An exhaustive search allows to compute $c_g(G)$ and a corresponding strategy for Alice in time $O(1)$. Roughly, the set of all games in $G$ can be described by one rooted tree with maximum degree at most $3q$ and depth $3q$. A classical dynamic-programming algorithm on this execution-tree can be used to compute the result in time $O(1)$.

	From now on, let us assume that $|V(G_2)|> 2q$. Note that, in this setting, as soon as Alice colours a vertex of $G_2$ (and she will always be able to do that in the strategies below because $|V(G_2)|> 2q$), all the red vertices of $K$ will belong to the same connected red component (since $G[K\cup R]= K\oplus R$). Moreover, in what follows, Alice will always colour at least $\left\lfloor \frac{|V(G_2)|}{2} \right\rfloor\geq q$ vertices in $G_2$, connected by a vertex of $K$, ensuring that the largest connected red component is always this one (the one containing all the red vertices of $G_2$) since $|V(G_1)|\leq q$.
	
	In what follows, we use the following variation of the  Maker-Breaker Largest Connected Subgraph game. Consider the following game that takes a graph $H$ and $X \subseteq V(H)$ as inputs. The game proceeds as the Maker-Breaker Largest Connected Subgraph game does, {\it i.e.}, Alice and Bob take turns colouring vertices of $G$ starting with Alice and with all the vertices being initially uncoloured. The difference lies in the objective of Alice. At the end of the game, the score achieved by Alice is the total number of red vertices that belong to the connected red components containing vertices of $X$. Intuitively, we see all the connected red components with at least one vertex in $X$ as a single connected red component. Let $c_g(H,X)$ be the largest integer $k$ such that Alice has a strategy to ensure a score of at least $k$ with input $(H,X)$, regardless of how Bob plays. Note that, by similar arguments as in the paragraph above, if $|V(H)|=O(1)$, then $c_g(H,X)$ (and a corresponding strategy for Alice) can be computed in time $O(1)$ for all $X \subseteq V(H)$.
	
    By the previous remark, $c_g(G_1,K)$ (and a corresponding strategy ${\cal S}_a^1$ for Alice) can be computed in time $O(1)$. By an exhaustive computation in constant time (since $|V(G_1)|=O(1)$), it is actually possible to consider all the strategies for Alice and Bob, including the ones where they may each skip at most one of their turns. If (in the variant game with input $(G_1,K)$) there exists a strategy for Alice guaranteeing her a score of at least $c_g(G_1,K)$, in which she skips one of her turns, and such that, if Bob skips a turn before Alice, then Alice can score at least $c_g(G_1,K)+1$ without skipping any of her turns, then let ${\cal S}_a^2$ be such a strategy for Alice. On the other hand, if (in the variant game with input $(G_1,K)$) there exists a strategy for Bob guaranteeing that Alice cannot score more than $c_g(G_1,K)$, in which he skips one of his turns, and such that, if Alice skips a turn before Bob, then Bob can guarantee that Alice scores at most $c_g(G_1,K)-1$ without skipping any of his turns, then let ${\cal S}^2_b$ be such a strategy for Bob. Note that, by definition, ${\cal S}_a^2$ and ${\cal S}^2_b$ cannot both exist simultaneously.
    
	Now, let us consider the following strategy ${\cal{S}}_b$ for Bob in $G$. Whenever Alice colours a vertex in $G_1$, Bob plays in $G_1$ following a strategy that ensures that Alice scores at most $c_g(G_1,K)$ in the variant game with input $(G_1,K)$. Whenever Alice colours a vertex in $G_2$, Bob colours any vertex of $G_2$ (if no such move is possible, Bob colours an arbitrary uncoloured vertex in $G$). This ensures that the largest connected red component is of order at most $c_g(G_1,K)+\left\lceil \frac{|V(G_2)|}{2} \right\rceil$. That is, $c_g(G) \leq c_g(G_1,K)+\left\lceil \frac{|V(G_2)|}{2} \right\rceil$.
	
	Let us also define the following strategy ${\cal{S}}_a$ for Alice in $G$. First, Alice colours the first vertex in $G_1$ that ensures her a score of at least $c_g(G_1,K)$ in the variant game with input $(G_1,K)$ (following strategy ${\cal S}_a^1$). Then, whenever Bob colours a vertex in $G_1$, Alice colours the vertex of $G_1$ following her strategy ${\cal S}_a^1$ to ensure a score $c_g(G_1,K)$ in the variant game with input $(G_1,K)$. Whenever Bob colours a vertex in $G_2$, Alice colours any vertex in $G_2$. If no such move is possible, Alice colours an arbitrary uncoloured vertex. This ensures that the largest connected red component is of order at least $c_g(G_1,K)+\left\lfloor \frac{|V(G_2)|}{2} \right\rfloor$ (recall that, since $|V(G_2)|\geq 2$, Alice colours at least one vertex in $G_2$). That is, $c_g(G) \geq c_g(G_1,K)+\left\lfloor \frac{|V(G_2)|}{2} \right\rfloor$.
	
	Note that the upper and lower bounds above match when $|V(G_2)|$ is even. Assume now that $|V(G_2)|$ is odd. We distinguish three cases. In all of the strategies below, the first player to colour a vertex in $G_2$ will colour at least $\left\lceil \frac{|V(G_2)|}{2} \right\rceil$ vertices in $G_2$.  
	\begin{itemize}
	    \item First, let us assume that the strategy ${\cal S}^2_a$ for Alice in $G_1$ defined above exists. We define Alice's strategy for $G$ as follows. Alice plays her first turns in $G_1$ following ${\cal S}^2_a$ until she can skip a turn in $G_1$ ({\it i.e.}, the first time she can skip a turn in $G_1$ while still guaranteeing a score of at least $c_g(G_1,K)$ in the variant game with input $(G_1,K)$). 
	    \begin{itemize}
	    \item If, in one of these rounds, Bob plays in $G_2$, then Alice first plays an extra turn in $G_1$ (following ${\cal S}^2_a$ that ensures her a score of at least $c_g(G_1,K)+1$ in the variant game with input $(G_1,K)$), and then, each time Bob plays in $G_1$, she plays in $G_1$ according to ${\cal S}^2_a$ in the variant game with input $(G_1,K)$, and each time Bob plays in $G_2$, she plays in $G_2$. 
	   
	    \item Otherwise, Bob also plays in $G_1$ until Alice can skip a turn in $G_1$. Then, once she can skip a turn in $G_1$ according to ${\cal S}^2_a$, Alice colours a vertex in $G_2$. From then, whenever Bob colours a vertex in $G_1$, she colours a vertex in $G_1$ following ${\cal S}^2_a$. Otherwise, she colours an arbitrary uncoloured vertex in $G_2$. 
	    \end{itemize}
	    In both cases, Alice scores at least $c_g(G_1,K)+\left\lceil \frac{|V(G_2)|}{2} \right\rceil$, matching the upper bound. 
	    
	    \item Second, let us assume that the strategy ${\cal S}^2_b$ for Bob in $G_1$ defined above exists. Note that ${\cal S}^2_a$ does not exist, so Alice cannot skip one turn in $G_1$ without decreasing her score in the variant game with input $(G_1,K)$. Bob plays his first turns in $G_1$ following ${\cal S}^2_b$ until he can skip a turn in $G_1$.
	    \begin{itemize}
	    \item If, in one of these rounds, Alice plays in $G_2$, then Bob first plays an extra turn in $G_1$ (following ${\cal S}^2_b$ that ensures him that Alice will score at most $c_g(G_1,K)-1$ in the variant game with input $(G_1,K)$). Then, whenever Alice plays in $G_1$, he continues to follow ${\cal S}^2_b$, and when Alice plays in $G_2$, Bob plays in $G_2$. 
	    \item Otherwise, Alice also plays in $G_1$ until Bob can skip a turn in $G_1$. Then, once he can skip a turn in $G_1$ according to ${\cal S}^2_b$, Bob colours a vertex in $G_2$. From then, whenever Alice colours a vertex in $G_1$, he colours a vertex in $G_1$ following ${\cal S}^2_b$. Otherwise, he colours an arbitrary uncoloured vertex in $G_2$.
	    \end{itemize}
	    In both cases, Alice scores at most 
	    $c_g(G_1,K)+\left\lfloor \frac{|V(G_2)|}{2} \right\rfloor$, matching the lower bound.
	    \item Finally, if none of the strategies ${\cal S}^2_a$ and ${\cal S}^2_b$ exist, the result depends on the parity of $|V(G_1)|$. Indeed, if Alice skips one turn in $G_1$, then Bob can ensure she scores at most $c_g(G_1,K)-1$ in the variant game with input $(G_1,K)$. On the other hand, if Bob skips one turn in $G_1$, Alice can score at least $c_g(G_1,K)+1$ in the variant game with input $(G_1,K)$. For Alice to ensure her upper bound and for Bob to ensure the lower bound, both of them will play in priority in $G_1$. That is, the first vertex of $G_2$ is coloured after all the vertices of $G_1$ have been coloured (and Alice has achieved a score of $c_g(G_1,K)$ in the variant game with input $(G_1,K)$). If $|V(G_1)|$ is even, Alice is the first player to colour a vertex in $G_2$, which allows her to score the upper bound $c_g(G_1,K)+\left\lceil \frac{|V(G_2)|}{2} \right\rceil$. Otherwise, Bob is the first player to colour a vertex in $G_2$, which implies that Alice can score at most the lower bound $c_g(G_1,K)+\left\lfloor \frac{|V(G_2)|}{2} \right\rfloor$.\qed 
	\end{itemize}
\end{enumerate}
\end{proof}


\section{Discussion and directions for further work}\label{section:conclusion}

A certain number of directions for further work on the Maker-Breaker Largest Connected Subgraph game seem particularly appealing to us. While some of the ones we mention are about tightening some of our results from the previous sections, others are original ones that are discussed only in this section.

\subsection{Differences between the two versions of the Largest Connected Subgraph game}

One direction for research could be to try to establish the significant differences between the Maker-Breaker Largest Connected Subgraph game and the Largest Connected Subgraph game. Some of our results in this work are already a step in that direction. For instance, in Section~\ref{section:complexity}, we showed  that the Maker-Breaker version remains PSPACE-complete when restricted to various classes of graphs, but we do not know whether the same holds for the Largest Connected Subgraph game in those classes of graphs. Lemma~\ref{lemma subgraph} draws another neat difference between the two versions of the game, as the outcome of the Largest Connected Subgraph game in a disconnected graph cannot be established as simply as in the Maker-Breaker version. This is because, in the latter version, Bob does not care about the structure induced by the blue vertices. However, in the Largest Connected Subgraph game, there are scenarios in which it is more favourable for Bob to play in a connected component $G_2$ different from the one $G_1$ that Alice just played in. This would be like skipping a turn in $G_1$, but playing an extra turn in $G_2$ (or playing first in $G_2$). Thus, to establish a result similar to Lemma~\ref{lemma subgraph} for the Largest Connected Subgraph game, one has to deal with the effects of skipping and playing extra turns, as well as Bob playing first, which seems like a tricky, yet interesting, aspect to study. 

\subsection{Types of strategies for the Maker-Breaker Largest Connected Subgraph game}

As we have seen in some graphs, notably in Section~\ref{section:bounds}, some optimal strategies for Alice ensure that the red subgraph is connected at all times. We believe it would be interesting to study a connected variant of the Maker-Breaker Largest Connected Subgraph game, in which Alice is always (except on her first turn) constrained to colour a neighbour of another red vertex, and the game ends when she cannot.
Consequently, we could define $c^c_g(G)$ as the maximum score Alice can achieve in $G$ when obliged to play in such a connected way. Clearly, $c^c_g(G) \leq c_g(G)$. We were able to observe that it is far from true that these two parameters are equal in general, even sometimes in quite simple graphs. As an illustration, the difference between both parameters is arbitrarily large for {\it king's grids} (strong products of two paths) with only two rows (denoted by $P_2 \boxtimes P_m$).

\begin{lemma}
For any $m \geq 1$, $c^c_g(P_2 \boxtimes P_m) = O(1)$ and $c_g(P_2 \boxtimes P_m)=m$. 
\end{lemma}

\begin{proof}
In the connected case, it is sufficient for Bob to colour the four vertices at distance $4$ from the first vertex coloured by Alice. In the non-connected case, each time Bob colours a vertex $v$, Alice colours the neighbour of $v$ in the other row (we say that $P_2 \boxtimes P_m$ has two rows and $m$ columns). At the end of the game, the red subgraph is connected, and so, Alice achieves a score of $m$.
\end{proof}

\subsection{Other classes of graphs}

Some of our results on particular classes of graphs leave open questions. Since the Maker-Breaker Largest Connected Subgraph game is PSPACE-complete in split graphs by Corollary~\ref{corollary:pspace-split}, and split graphs have diameter at most~$3$, there is the question of whether it is hard to compute $c_g$ for graphs of diameter~$2$. Regarding the results from Section~\ref{section:sparse}, recall that~\cite{paper1} provides a linear-time algorithm for the Largest Connected Subgraph game in cographs. One question is whether this result can be extended to $(q,q-4)$-graphs, as we did in Theorem~\ref{theorem:qgraphs} for the Maker-Breaker version.

Regarding determining $c_g$ for other graph classes, an appealing direction could be to consider standard graph classes such as trees. From~\cite{paper1}, we have that $c_g(P_n)=2$ for any path $P_n$ of order $n\geq 3$, and we believe that understanding the game in larger subclasses of trees such as caterpillars and subdivided stars is not so difficult, but requires a lot of work to prove, for a not so substantial result. Thus, we think it would be most interesting to study the class of trees rather than its subclasses. Other natural graph classes to be investigated are graph products. For instance, we wonder whether $c_g(Q_n)$ can be easily determined for a hypercube $Q_n$ (where, recall, $Q_2$ is the cycle $C_4$ of length~$4$,
and, for every $n > 2$, the hypercube $Q_n$ is the Cartesian product $Q_{n-1} \square P_2$ of $Q_{n-1}$ and the path $P_2$ of order~$2$).
We also wonder about different types of grids, which, in the Largest Connected Subgraph game, seem hard to comprehend~\cite{paper1}. Another point for considering such graphs is that grids are natural structures to play on in several types of games, as illustrated by Hex. To give some insight into what can be done in the Maker-Breaker version, we finish off with some partial results on grids in the rest of this section. We first consider {\it hexagonal grids}.

\begin{figure}[t!]
\begin{center}
\includegraphics[width=5cm]{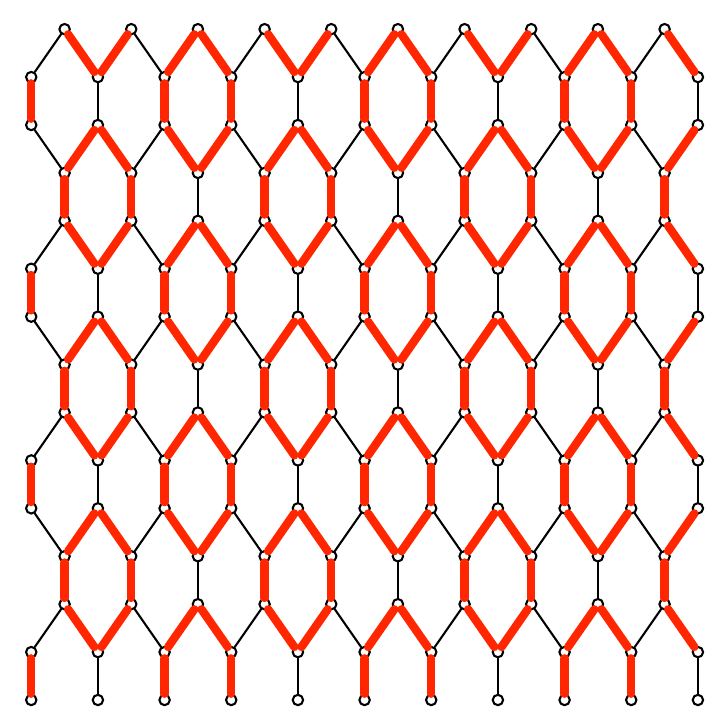}
\end{center}
\caption{Illustration of the infinite Hexagonal grid $H_{\infty}$ in the proof of Proposition~\ref{prop:hexGrid}. The connected red subgraphs are vertex-disjoint $6$-cycles covering the vertices of $H_{\infty}$. The black edges induce a matching of $H_{\infty}$.}
\label{fig:hexGrid}
\end{figure}

\begin{proposition}\label{prop:hexGrid}
If $G$ is a finite subgraph of the infinite hexagonal grid, then $c_g(G) \leq 6$. 
\end{proposition}

\begin{proof}
Let $H_{\infty}$ be the infinite hexagonal grid as partially shown in Fig.~\ref{fig:hexGrid}. Let $(C_i)_{i \in \mathbb{N}}$ be the set of vertex-disjoint subgraphs of $H_{\infty}$ depicted in red in Fig.~\ref{fig:hexGrid}. 
For any $i \in \mathbb{N}$, $C_i$ induces a cycle of order $6$ and $(V(C_i))_{i \in \mathbb{N}}$ is a partition of $V(H_{\infty})$. Furthermore, $M=E(H_{\infty})\setminus (\bigcup_{i \in \mathbb{N}} E(C_i))$ (black edges in Fig.~\ref{fig:hexGrid}) is a matching of $H_{\infty}$. Note also that, for any $i \neq j$, every path from a vertex of $C_i$ to a vertex of $C_j$ contains an edge in $M$ (since, for every subgraph $C_i$, the edges adjacent to a vertex of $V(C_i)$, but not in $E(C_i)$, are by definition in $M$). 

Let $G$ be any finite subgraph of $H_{\infty}$, and let $M'$ and $(V(C'_i))_{i \in \mathbb{N}}$ be the restrictions of $M$ and $(V(C_i))_{i \in \mathbb{N}}$ to $G$, respectively. Consider the following strategy for Bob in $G$. First, note that, for any vertex $v\in V(G)$, there is at most one edge $uv\in M'$ incident to $v$ since $M$ is a matching of $H_{\infty}$. Thus, each time Alice colours a vertex $v$, Bob colours the vertex $u$ such that $uv\in M'$, if it exists and it is uncoloured, and if not, then he colours an arbitrary uncoloured vertex in $G$. Let us show that Bob's strategy ensures that Alice cannot create a connected red component of order more than $6$. Towards a contradiction, let us assume that Alice creates a connected red component $S$ of order at least $7$. Then, there exist $u,v \in S$ and $i \neq j$ such that $u \in V(C'_i)$ and $v \in V(C'_j)$ (because the $C_k$'s partition the vertex-set of $H_{\infty}$ and each $C_k$ has order $6$). As mentioned above, every path between $u$ and $v$ contains an edge of $M'$, and so, by Bob's strategy, a vertex of this path was coloured by Bob, contradicting that $u$ and $v$ belong to the same connected red component.
\end{proof}

Through a tedious case analysis, it might be possible to prove that Proposition~\ref{prop:hexGrid} is sharp, that is, that there exists a finite subgraph $G$ of $H_{\infty}$ such that $c_g(G) = 6$. We also consider {\it Cartesian grids}, which are the Cartesian product $P_n \square P_m$ of $P_n$ and $P_m$. For these grids, we provide the following upper bound.

\begin{proposition}\label{prop:square-grids}
For $n \leq m$, $c_g(P_n \square P_m) \leq 2n$. 
\end{proposition}

\begin{proof}
Consider an $n \times m$ grid $P_n \square P_m$ with $n$ rows and $m$ columns (with left, right, higher, and lower being defined naturally). Consider the following strategy for Bob. When Alice colours a vertex $v$, if the right neighbour $u$ of $v$ exists and is uncoloured, then Bob colours $u$, otherwise, Bob colours the left neighbour of $v$ if it exists and is uncoloured, and otherwise, Bob colours an arbitrary uncoloured vertex. 

The above strategy for Bob is well-defined and ensures that no three consecutive vertices in a row are ever red (see the case of paths in~\cite{paper1} for more details). This ensures that, for any strategy of Alice, any connected red component has at most $2$ vertices in each row, hence, proving the proposition. Indeed, consider a largest connected red component $S$ at the end of the game. Towards a contradiction, assume that there exists a row whose intersection with $S$ contains strictly more than $2$ vertices. Then, the restriction of $S$ to this row $r$ induces at least two connected red components since there cannot be three consecutive red vertices in a same row.

Let $x,y \in V(P_n \square P_m)$ be two vertices of $S$ in a same row $r$ such that they are in different components of the intersection of $r$ and $S$. Let $P$ be any shortest red path from $x$ to $y$ (it exists since $S$ is connected). By the definitions of $x$ and $y$, $P$ contains a vertex in a row above or below $r$. W.l.o.g., suppose $P$ contains a vertex in a row above $r$. Let $w$ be a highest vertex of $P$, {\it i.e.}, no vertex of $P$ is in a higher row than $w$. Let $r'$ be the row of $w$. Then, $r'$ is higher than $r$. Thus, $w\notin \{x,y\}$, and so, $w$ has degree two in $P$. Then, since $w$ is the highest vertex of $P$, there exists a red neighbour $z$ of $w$ in $V(P) \cap r'$. W.l.o.g., say that $z$ is to the right of $w$. Then, $z\notin \{x,y\}$ either, and so, $z$ has degree two in $P$. Note that there can be no vertex of $P$ directly to the right of $z$ nor directly to the left of $w$, since otherwise, there would be $3$ consecutive red vertices. Furthermore, there can be no vertex of $P$ above $w$ or $z$ since $w$ is a highest vertex of $P$. Hence, since $w$ and $z$ have degree two in $P$, the vertex directly below $w$, and the one directly below $z$, must also be in $P$. Then, the fact that $w$ and $z$ belong to $P$ contradicts that $P$ is a shortest path.
\end{proof}

Regarding Proposition~\ref{prop:square-grids}, we would be interested in knowing the precise value of $c_g(P_n \square P_m)$ in general. One issue we ran into is the fact that Alice can play in a non-connected way (recall the notion of connected moves discussed earlier), and it is not clear how Bob should anticipate to prevent connected red components to merge later on.
Let us mention, however, that if Alice plays in a connected way in a Cartesian or king's grid, then the game becomes quite similar to the Angel and Devil Problem of Conway~\cite{winningways}. Optimal strategies for the devil in that game~\cite{winningways} allow to prove that $c_g^c(P_n \square P_m)$ is bounded above by an absolute constant.

\bibliography{LCSG}
\bibliographystyle{plain} 

\end{document}